\documentclass[3p]{elsarticle}

\usepackage{moreverb}
\usepackage{array}
\newcolumntype{L}[1]{>{\raggedright\let\newline\\\arraybackslash\hspace{0pt}}m{#1}}
\newcolumntype{C}[1]{>{\centering\let\newline\\\arraybackslash\hspace{0pt}}m{#1}}
\newcolumntype{R}[1]{>{\raggedleft\let\newline\\\arraybackslash\hspace{0pt}}m{#1}}

\usepackage[hyphens]{url}

\newcommand\BibTeX{{\rmfamily B\kern-.05em \textsc{i\kern-.025em b}\kern-.08em
T\kern-.1667em\lower.7ex\hbox{E}\kern-.125emX}}

\usepackage{amsfonts,amsmath,mathrsfs,amssymb,latexsym,epsfig}

\usepackage{graphicx,psfrag,pstricks,verbatim,wasysym}
\usepackage{algorithm}
\usepackage{bclogo}
\usepackage[noend]{algpseudocode}

\newcommand{\I}{\textit{I}}

\newcommand{\II}{\textit{II}}
\newcommand{\III}{\textit{III}}

\newcommand{\dsum}{\textstyle \sum}
\newcommand{ \bs}[1]{\boldsymbol{#1}}

\newcommand{\PP}{\hbox{I\kern-.2em\hbox{P}}}
\newcommand{\HH}{\hbox{I\kern-.2em\hbox{H}}}

\newcommand{\RR}{\hbox{I\kern-.2em\hbox{R}}}

\usepackage[latin1]{inputenc}
\usepackage{commath}
\usepackage{subfig}
\usepackage{epsfig}
\usepackage{booktabs}
\usepackage{tabularx}
\usepackage{listings}
\usepackage{mathtools}
\usepackage{multirow}
\usepackage{tabu}
\usepackage{romanbar}
\usepackage{xfrac}
\usepackage{color}

\usepackage{subfig}
\usepackage{float}
\usepackage{placeins}
\usepackage{rotating}
\usepackage{pdflscape}

\usepackage{etoolbox}
\apptocmd{\sloppy}{\hbadness 10000\relax}{}{}

\newtheorem{thm}{Theorem}
\newtheorem{lem}[thm]{Lemma}
\newdefinition{rmk}{Remark} \newproof{pf}{Proof}

\newcommand{\proofofref}{}
\newproof{zproofof}{Proof of Lemma \proofofref}
\newenvironment{pol}[1]
 {\renewcommand{\proofofref}{#1}\zproofof}
 {\endzproofof}

\newproof{tproofof}{Proof of Theorem \proofofref}
\newenvironment{pot}[1]
 {\renewcommand{\proofofref}{#1}\tproofof}
 {\endzproofof}

\usepackage{setspace}
\begin{document}
\title{An Efficient and Accurate Method for Modeling Nonlinear Fractional Viscoelastic Biomaterials}

\author[umich]{Will Zhang}
\author[kcl]{Adela Capilnasiu}
\author[graz]{Gerhard Sommer}
\author[graz,ntnu]{Gerhard A. Holzapfel}
\author[umich,umich2,kcl]{David A. Nordsletten}

\address[umich]{Department of Biomedical Engineering, University of Michigan, Ann Arbor, USA}
\address[umich2]{Department of Cardiac Surgery, University of Michigan, Ann Arbor, USA}
\address[graz]{Institute of Biomechanics, Graz University of Technology, AT}
\address[ntnu]{Department of Structural Engineering, Norwegian University of Science and Technology, Trondheim, NO}
\address[kcl]{School of Biomedical Engineering and Imaging Sciences, King's College London, London, UK}

\begin{abstract}
Computational biomechanics plays an important role in biomedical engineering: using modeling to understand pathophysiology, treatment and device design.
While experimental evidence indicates that the mechanical response of most tissues is viscoelasticity, current biomechanical models in the computation community often assume only hyperelasticity.
Fractional viscoelastic constitutive models have been successfully used in literature to capture the material response. However, the translation of these models into computational platforms remains limited. 
Many experimentally derived viscoelastic constitutive models are not suitable for three-dimensional simulations.
Furthermore, the use of fractional derivatives can be computationally prohibitive, with a number of current numerical approximations having a computational cost that is $ \mathcal{O} ( N_T^2) $ and a storage cost that is $ \mathcal{O}(N_T) $ ($ N_T $ denotes the number of time steps).
In this paper, we present a novel numerical approximation to the Caputo derivative which exploits a recurrence relation similar to those used to discretize classic temporal derivatives, giving a computational cost that is $ \mathcal{O} (N) $ and a storage cost that is fixed over time.
The approximation is optimized for numerical applications, and the error estimate is presented to demonstrate efficacy of the method.
The method is shown to be unconditionally stable in the linear viscoelastic case. 
It was then integrated into a computational biomechanical framework, with several numerical examples verifying accuracy and computational efficiency of the method, including in an analytic test, in an analytic fractional differential equation, as well as in a computational biomechanical model problem.
\end{abstract}

\begin{keyword}
Caputo derivative; viscoelasticity; solid mechanics; computational biomechanics; large deformation
\end{keyword}

\maketitle


\section{Introduction}

    Biomechanical analysis has become an area of increasing emphasis in biomedical engineering applications.
    It has been used for understanding a diverse range of complex problems, such as the surgical intervention strategies for arterial bifurcation \cite{moore2010coronary}, the impact of left ventricular assist devices \cite{mccormick2011modeling,mccormick2014computational}, the risk of aortic aneurysm rupture~\cite{gasser2010biomechanical,gasser2014novel}, the repair of the mitral valve in the heart~\cite{rausch2013mechanics,dasi2017pursuit}, as well as the surgical aids for brain surgeries \cite{miller2011biomechanical}.
    Through biomechanics simulation, devices and procedures can be optimized and personalized to improve design, therapy and patient outcomes~\cite{chabiniok2016multiphysics}.
    A key component to these types of biomechanical assessments is identifying appropriate constitutive equations that define the material response to physiological and extra-physiological loading.
    While hyperelastic material models have been developed to describe a range of tissues~\cite{fung1993biomechanics,holzapfel2006mechanics}, viscoelasticity has been acknowledged as an important component of mechanical response in brain matter \cite{streitberger2010vivo,li2016elastic,sack2009impact}, connective tissues \cite{dunn1983viscoelastic,lanir1983constitutive}, fracturing of bone tissues \cite{wu2011viscoelastic}, as well as many cardiovascular soft tissues \cite{wang2016viscoelastic}, including arterial tissues \cite{holzapfel2002structural} and myocardium \cite{sommer2015biomechanical,sommer2015quantification}. 
    Despite viscoelasticity in soft tissues being widely studied experimentally, its use in computational biomechanics simulations remains relatively uncommon, likely due, in part, to the challenges to its implementation.
    
    Viscoelasticity and viscoelastic theory has been a topic of interest since the 19th century~\cite{maxwell1867iv,wiechert1893gesetze}.
    Early forms of viscoelastic constitutive models used spring and dashpot~\cite{wilhelm1975viscoelasticity} constructs, e.g, the Maxwell~\cite{christensen1980nonlinear} and the Kelvin-Voigt models~\cite{meyers2014mechanical}.
    To better explain experimental observations, authors began to combine these simple units in series and/or parallel to form more complex systems, such as the generalized Maxwell~\cite{wiechert1893gesetze} or generalized Kelvin-Voigt models \cite{casula1992generalized}.
    These early theoretical works were successfully used to explain linear viscoelastic behavior a number of materials, including muscle tissues~\cite{levin1927viscous}.
    A more generalized one dimensional theory was proposed by Green and Tobolsky~\cite{green1946new} to explain relaxing media.
    Extension of these early concepts of viscoelasticity to nonlinear continuum mechanics was addressed by a number of analytic mechanicians, including Truesdell~\cite{truesdell1955simplest}, Green and Rivlin~\cite{green1957mechanics,green1959mechanics} and Bernstein, Kearsley, and Zapas~\cite{bernstein1963study} (see~\cite{coleman1961foundations} for a review on early works investigating viscoelasticity).
    Nonlinear viscoelastic theory has been applied to modelling of elastomers~\cite{christensen1980nonlinear}, rubber materials~\cite{lubliner1985model}, fiber reinforced materials~\cite{holzapfel2001viscoelastic},  myocardium~\cite{gueltekin2016orthotropic} as well as brain tissue~\cite{budday2017rheological,budday2017viscoelastic}.
    Extension of these approaches to computational mechanics has been in the works of Simo~\cite{simo1987fully}, Le Tallec~\cite{le1993three}, Holzapfel~\cite{holzapfel1996new,holzapfel2001viscoelastic,holzapfel2002structural,gueltekin2016orthotropic}, and De Buhan and Frey~\cite{de2011generalized} amongst others.
    
    While much of the early developments in viscoelasticity and its numerical application focused on analogues to spring-dashpot systems, fractional viscoelastic behavior can be seen increasingly in modeling tissue properties~\cite{magin2006fractional,holm2019waves,freed2006fractional}.
    Use of these models has been recently presented in lung tissue~\cite{birzle2019viscoelastic}.
    The advantage of fractional viscoelasticity is that the method provides a spectrum of relaxation times parameterized through the chosen fractional derivative order, $ \alpha $.
    A simplification of the generalized Maxwell approach with infinite Maxwell elements, the fractional approach has the advantage of being compact, requiring few parameters and being effective at characterizing hierarchical viscoelastic materials exhibiting power-law behavior~\cite{magin2004fractional,hilfer2000fractional,sverregrimnes2014bioimpedance,lakes1999viscoelastic}.
    These benefits are challenged by the need to numerically approximate the fractional derivative, which is often done using cumulative integral approximations that require the storage of all past values of stress and increasingly costly computations in time.
    Development of efficient approaches to approximating fractional derivatives by Yuan and Agrawal~\cite{yuan2002numerical},  Diethelm~\cite{diethelm2008investigation} and Birk and Song~\cite{birk2010improved} have enabled the efficient extension of fractional models to solid mechanics~\cite{zopf2015comparison}.
    However, how these approximations perform in a computational biomechanics context using finite elements presents additional challenges and questions.

    In this work, we present a numerical formulation for the constitutive modeling of viscoelastic materials based on a novel and efficient recursive approximation to evaluate fractional viscoelasticity. 
    The proposed approach, utilizes a Prony series (or generalized Maxwell model) as in~\cite{diethelm2009improvement,birk2010improved,zopf2015comparison} with an additional dashpot element and optimization process to making the computational and storage costs consistent with standard transient methods.
    An error estimate is derived for the class of Prony approximations, detailing the improved computational accuracy of the proposed method (as well as those previously proposed).
    The numerical method is then integrated into our finite element framework~\cite{lee2016multiphysics} and shown to be unconditionally stable in the linear viscoelastic setting.
    Results from this analysis are demonstrated through a series of numerical examples, demonstrating appropriateness of the developed bounds, accuracy and numerical efficiency.
    While demonstrated in a specific biomechanics setting, this generalized approach lends itself to a wide variety of modeling problems in soft tissues as well as other fractional-time derivative problems.
    
    The paper reviews current techniques for approximating the Caputo derivative, introduces our approximation as well as its error estimates (section~\ref{meth:frac_derivatives}). 
    In section~\ref{meth:solid_mech}, the paper reviews mechanical formulations for viscoelasticity, including fractional approaches, for nonlinear mechanics.
    The introduced models and approximations are then integrated into a nonlinear mechanics finite element framework (section~\ref{meth:solid_mech_num}) which is shown to be unconditionally stable in the linear viscoelastic case (see Appendix~\ref{app:stability}). 
    Convergence, accuracy and computational costs are then evaluated in an analytic example (section~\ref{sect:poly_example}), a simplified fractional diffusion system (section~\ref{sect:gao_example}) as well as in a biomechanical example of liver tissue (section~\ref{sect:liver_example}). 
    The paper finishes with a discussion of key results and conclusions.


\section{Fractional Derivatives and their Approximations} \label{meth:frac_derivatives}

\subsection{Caputo fractional derivative}

The concept of fractional calculus started with questions about about the generalization of integral and differential operators by L'Hospital and Leibniz~\cite{ross1977development} from the set of integers to the set of real numbers.
Subsequently, many prominent mathematicians focused on fractional calculus (for reviews, see, Ross~\cite{ross1977development} and Machado~\cite{machado2011recent}).
Within the field, many different definitions for fractional differential and integral operators of arbitrary order have been introduced~\cite{podlubny1998fractional}.
In this work, we focus on the Caputo definition where the fractional derivative, $ D_t^\alpha $ (with $ \alpha > 0 $), of a $ n $-times differentiable function $ f  $ can be written as,
\begin{equation} \label{eq_caputo_def}
  D_t^\alpha f = \frac{1}{\Gamma( \lceil \alpha \rceil - \alpha)} \int_0^t \frac{ f^{\lceil \alpha \rceil} (s)}{(t-s)^{1+\alpha - \lceil \alpha \rceil}} \dif s,  
\end{equation}
where $ \lceil \alpha \rceil $ denotes the ceiling of $ \alpha $.
The fractional derivative (for $ \alpha \notin \mathbb{N}$), unlike integer derivatives, is a convolution of the past behavior of the function, making it particularly useful for constitutive modeling of viscoelastic materials.
The Caputo form (which can be related to other formulations like the Riemann-Louiville fractional derivative) has the advantage that
\begin{equation} \label{eq_caputo_def2}
  D_t^\alpha c = 0, \quad \text{ for any constant } c \in \mathbb{R}.
\end{equation}
For later analysis, we will instead assume that the functions $ f $ are continuous and $ n+1 $-times differentiable~\cite{diethelm2005algorithms}, which allows the use integration by parts to rewrite the Caputo derivative as
\begin{equation}
    D_t^\alpha f  = \frac{1}{\Gamma(1+n-\alpha)} \left (  
        t^{n-\alpha} f^n (0) + \int_0^t (t-s)^{n-\alpha} f^{n+1} (s) \dif s
    \right).
\end{equation}
This form eliminates the weak singularity in the integrand of Eq.~\eqref{eq_caputo_def}, which is convenient in our analysis.

\subsection{Current methods for numerical approximation of the Caputo derivative}\label{sect:current_frac_meths}

Fractional derivative operators and methods for approximating their effects are well studied in literature. 
For some recent review of this subject, readers may consider the works of Zeid \cite{zeid2019approximation}, Guo \cite{guo2015numerical}, Weilbeer \cite{weilbeer2005efficient} and Podlubny \cite{podlubny1998fractional}. 
Approaches in the literature can be categorized by their method of integration, either cumulative  or recursive.
Cumulative approaches approximate the Caputo derivative by directly numerically approximating the integrals in Eq.~\eqref{eq_caputo_def} or~\eqref{eq_caputo_def2} using some form of a weighted sum over a set of discrete time points.  
In contrast, recursive approaches first approximate the Caputo derivative by transforming it into a series of ordinary differential equations that have standard discretizations.

Suppose we seek to approximate the $ \alpha^{th}$-order fractional derivative, $ D_t^\alpha f $ ($0 < \alpha < 1 $), of a function $ f : [0,T] \to \mathbb{R} $ at a series of $N_T$ times points 
$$
    0=t_0 < t_1 < \cdots < t_{N_T} = T, \quad \text{ where } t_n = n \Delta_t,
$$    
$ \Delta_t$ denotes the time step size, and that $ f_n = f(t_n) $ is easily computed given some prior information from previous time points. 
Considering cumulative approaches, a straightforward method to computing Caputo derivative is through a Riemann sum evaluated using the midpoint point (MP) rule, i.e.
\begin{equation}\label{eqn:basiccaputo}
   \text{D}_{n,\text{MP}}^\alpha f = \frac{1}{\Gamma(1-\alpha)} \sum_{i=1}^n \frac{1}{\left(t-(i-\frac{1}{2})\Delta_t\right)^\alpha} 
    (f_{k} - f_{k-1}).
\end{equation}
Here the integral in Eq.~\eqref{eq_caputo_def} is approximated at the midpoint of each time interval with the derivative computed by central difference.
The weights of integration decays temporally.
Other cumulative methods differ by the choice of weights to improve accuracy or achieve higher order convergence. 

Another common method for constructing the fractional order derivative is using the Gr\"unwald-Letnikov (GL) derivative~\cite{Scherer2011}, which extends on the standard finite difference definition for the derivative of functions.  
Using the definition from Scherer \emph{et al.}~\cite{Scherer2011}, the GL derivative with corrections to match the Caputo derivative is
\begin{equation}\label{eq_gl}
    \text{D}_{n,\text{GL}}^\alpha f =
            \frac{1}{\Delta_t^\alpha} \left [ 
                f_{n} - \sum_{m=1}^n c_m^\alpha f_{n-m}
            \right]
            -  \frac{t_n^{-\alpha}}{\Gamma(1-\alpha)} f_0,
\end{equation}
where $ C_m^\alpha = C_{m-1}^\alpha (1 - [\alpha+1]/ m) $, $ C_1^\alpha = \alpha $, is the recursive definition of the real numbered binomial coefficient in the finite difference formula.
Another example is given by Diethelm \textit{et al.} \cite{diethelm1997generalized, diethelm2005algorithms} where the trapezoidal rule was employed, i.e.
\begin{equation} \label{eqn:diethelgrunwald}
\begin{aligned}
    \text{D}_{n,D}^\alpha f =& \left[\frac{\Delta_t^{-\alpha}}{\Gamma(2-\alpha)}\right] \sum_{m=0}^n a_{m,n} \left(f_{n-m} - \sum_{k=0}^{\lceil\alpha\rceil}\left[\frac{(n-m)^k \Delta_t^k}{k!}\right]f_0^{(k)}\right),\\
    a_{m,n} =&
    \begin{cases}
1, &\text{if} \ m = 0,\\
(m+1)^{1-\alpha} - 2m^{1-\alpha}+(m-1)^{1-\alpha}, &\text{if} \ 0<m<n,\\
(1-\alpha)n^{-\alpha} - n^{1-\alpha}+(n-1)^{1-\alpha}, &\text{if} \ m=n,
\end{cases}
\end{aligned}
\end{equation}
    where $a_{m,n}$ are the weights for the past values of $f_n$, and $f_0^{(k)}$ is the $k$-th derivative at 0. Richardson interpolation can also be used by selecting a subset for different size of $\Delta_t$ to improve accuracy. 
    Gao \textit{et al. } \cite{gao2012finite} developed a similar method for solving fractional diffusion equations, \emph{i.e.} 
\begin{equation} \label{eqn:gaoderivative}
\begin{aligned}
    &\text{D}_{n,\text{Gao}}^\alpha f = \frac{\Delta_t^{-\alpha}}{\Gamma(2-\alpha)}\left[a_0^\alpha f_n - \sum_{k=1}^{n-1}\left( a_{k-i-1}^\alpha - a_{k-i}^\alpha \right)f_i - a_{k-1}^\alpha f_0\right], \\
    &\quad a_i^\alpha = (i+1)^{1-\alpha} - i^{1-\alpha}, \quad i \geq 1.
\end{aligned}
\end{equation}
    Other cumulative methods have been presented by Murio \textit{et al.} \cite{murio2008implicit}, Yang \textit{et al.} \cite{yang2010numerical} as well as Liu et al. \cite{liu2007stability}, though many more have been published in the literature.
    
    All the cumulative methods presented here require significant storage costs and have significant computational costs, especially for complex three dimensional tissue simulations.
    Storage/memory costs are typically $ \mathcal{O}(N_T) $ and the computational complexity is often $ \mathcal{O}(N_T^2) $, or by best reports $ \mathcal{O}(N_T \log(N_T))$, due to the need to sum over the entire time domain.
    In contrasting, the computational cost of transient elasticity problems is independent of the number of time steps.  
    
    One strategy to cope with this limitation is to exploit the fact that the integration weights for numerical approximation of fractional derivatives exhibit a 'fading memory', where the weights tend to decay with the number of time steps (or length of time). 
    Hence, some have proposed truncating the approximation, including only the last few weights.
    However, the accuracy of this approach is fundamentally dependent on the behavior of $ f $ and the fractional derivative order $ \alpha $. 
    For example, the weights decay at slower rates for lower $ \alpha $ values, requiring more and more terms to avoid inaccuracies due to truncation.
    More difficult is the dependence on the function $ f $, which is not known \emph{a priori} in computational simulations.
    A classic example is to consider a stress relaxation experiment, where an acute action occurring at the beginning of the experiment dictates the response over the entire time course.
    
    Another strategy to reduce computational costs is by use of recursive methods.
    Such methods mostly follow a similar approach~\cite{zopf2015comparison}.
    Namely, the Caputo derivative is first approximated using decaying exponentials (or, equivalently, a Prony-series). 
    This approximation step is posited differently depending on the authors. However, a natural rationale for this approximation stems from the relationship between the relaxation spectrum of the Caputo derivative and decaying exponentials via the Laplace transform, i.e.
    \begin{equation}\label{eq_frac_trans}
        \frac{(t-s)^{-\alpha }}{\Gamma(1-\alpha)} 
        =
        \int_0^\infty \left[\frac{\sin(\pi \alpha) z^{\alpha-1} }{\pi}  \right] \exp[(s-t)z] \dif z.
    \end{equation}
    The integral above can be approximated by summation at a set of points $z$ with the associated quadrature weights (in the bracketed term).
    Inserting into Eq.~\eqref{eq_caputo_def}, we arrive at an equivalent form for a Prony-series comprised of decaying exponentials that can be recast as a series of ordinary differential equations (shown in the following section).
    One of the earlier works presenting a recursive method for approximating the Caputo derivative was by Yuan-Agarwal~\cite{yuan2002numerical}.
    In this case, authors recast the integral based on the above and used Gauss-Laguerre quadrature for integration.
    Though they presented the idea in the context of a specific fractional differential equation, the idea was clearly generalizable.
    

    The Yuan-Agrawal method has been critized for slow convergence for certain values of $\alpha$. Investigation by Diethelm \cite{diethelm2008investigation} found the integrand of the the Yuan-Agrawal integrand method to be non-smooth at 0 for $\alpha \neq 1/2 + n$, contrary to the assumption made in the Gauss-Laguerre quadrature. This results in poor convergence for $\alpha\to0$ and $\alpha\to\infty$.
    Instead, Diethelm suggested is his method \cite{diethelm2009improvement} that Gauss-Laguerre be replaced by Gauss Jacobi integration by mapping the range of integration from $  [0, \infty)$ to $ [-1,1]$ through a change of variables.  
    Aiming for further improvements in the accuracy to these methods, Birk and Song \cite{birk2010improved} further extended the Diethelm method \cite{diethelm2009improvement} by altering the mapping  onto the $ [-1,1]$ interval and considering the Fourier transform of the fractional operator, showing generally comparable or superior results.
    
    Although the quadrature scheme of these approaches covers the full range of frequency, this is not necessarily an advantage. The range of relevant frequencies can be limited based on the corresponding problem and by the size of $\Delta_t$ used for simulation. 
    Moreover, as $ \alpha \to 1 $, the approximation must tend to approach the first order derivative, which becomes increasingly challenging to approximate.
    In the following sections, we expand on this Prony-based approach by formulating the error estimates, generalizing its form and proposing a new approach for optimizing parameters of the weights and time constants of the associated Prony series.
    


\subsection{Prony approximation of the Caputo derivative}
\label{meth:frac_approx}
%
    The Prony series approach has been used a number of times in literature~\cite{yuan2002numerical,birk2010improved,peter2013generalized, potts2010parameter,diethelm2008investigation,diethelm2009improvement}.
    Suppose that the fractional derivative can by approximated using $ N $ Maxwell elements in parallel, each with their own weight $ \beta_k $ and time constant $ \tau_k $ ($k=1 \ldots N $).  
    Then a Prony-based approximation $ \hat{D}_t^{\alpha} f $ to $ D_t^\alpha f$ can be written as,
\begin{equation}\label{eq_Prony_approx}
   \hat{D}_t^{\alpha} f 
      := 
      \beta_0 f^\prime(t)  
      + 
      \sum_{k=1}^N \int_0^t \beta_k \exp \left[ \frac{s-t}{\tau_k} \right] f^\prime (s) \dif s,
\end{equation}
or equivalently as
\begin{equation}\label{eq_Prony_approx_series}
   \hat{D}_t^{\alpha} f 
      := 
      \beta_0 f^\prime(t)  + \sum_{k=1}^N q_k (t), 
      \quad 
      q_k^\prime(t) + \frac{1}{\tau_k} q_k(t) = \beta_k f^\prime(t),
\end{equation}
where we introduce $ N $ intermediate variables $ q_k $. 
While Eq.~\eqref{eq_Prony_approx} requires further approximation of the integral over the time domain, Eq.~\eqref{eq_Prony_approx_series} enables discretization of the intermediate variables $ q_k $.
Eq.~\eqref{eq_Prony_approx_series} requires $2N+1 $ parameters, $  \bs{\theta}_\alpha = \{\tau_1, \ldots, \tau_N, \beta_0, \beta_1,\ldots, \beta_N \} $.
The approximation error of equations~\eqref{eq_Prony_approx} and~\eqref{eq_Prony_approx_series} is shown in Lemma~\ref{lem_Prony_error}:
\begin{lem}\label{lem_Prony_error}
Let $ f: [0,T] \to \mathbb{R} $ be a real function on the domain [0,T]. If $ f^\prime(0) < \infty $ and $ f^{\prime \prime} \in L^1(0,T) $, then
\begin{equation*}
    (D_t^\alpha f - \hat{D}_t^\alpha f) = \varepsilon (t) f^\prime(0) + \int_0^t \varepsilon(z) f^{\prime \prime} (t-z) \dif z,
\end{equation*}
and for any $ t \in [0,T] $
\begin{equation*}
    \left | D_t^\alpha f - \hat{D}_t^\alpha f \right |
    \le 
    \| \epsilon \|_{0,\infty}
    \left [ 
    | f^\prime (0) | +
    \| f^{\prime \prime} \|_{0,1}
    \right ],
\end{equation*}
where the truncation error $ \varepsilon $ is defined as
\begin{equation} \label{eq_Prony_error}
    \varepsilon(z) := \frac{z^{1-\alpha}}{\Gamma(2-\alpha)} - \beta_0 + \sum_{k=1}^N \beta_k \tau_k ( \exp(-z / \tau_k ) - 1).
\end{equation}
\end{lem}

\begin{pol}{\ref{lem_Prony_error}}
Applying integration by parts to Eq.~\eqref{eq_Prony_approx}, and the fundamental theorem of calculus to the leading first order derivative term,
the Prony approximation to the fractional derivative can be re-written as,
\begin{equation*}
    \hat{D}_t^\alpha f  = \kappa(0,t) f^\prime (0) 
       + 
       \int_0^t \kappa(s,t) f^{\prime \prime} (s) ds,
\end{equation*}
where  $  \kappa(s,t) = \beta_0 - \sum_{k=1}^N \beta_k \tau_k ( \exp([s-t]/\tau_k) - 1 ) $.
Subtracting this from Eq.~\eqref{eq_caputo_def2} lead to
\begin{eqnarray}
    (D_t^\alpha f - \hat{D}_t^\alpha f)
    & = & 
    \frac{1}{\Gamma(2-\alpha)} \left (  
        t^{1-\alpha} f^\prime (0) + \int_0^t (t-s)^{1-\alpha} f^{\prime \prime} (s) \dif s
    \right)
    -
    \kappa(0,t) f^\prime (0) 
    - 
    \int_0^t \kappa(s,t) f^{\prime \prime} (s) \dif s
    \nonumber \\
    & = &
    \left( \frac{t^{1-\alpha}}{\Gamma(2-\alpha)} - \kappa(0,t) \right) f^\prime (0)
    +
    \int_0^t \left ( \frac{(t-s)^{1-\alpha}}{\Gamma(2-\alpha)} - \kappa(s,t) \right) f^{\prime \prime} (s) \dif s
    \nonumber \\
    & = &
    \left( \frac{t^{1-\alpha}}{\Gamma(2-\alpha)} - \kappa(0,t) \right) f^\prime (0)
    +
    \int_0^t \left ( \frac{z^{1-\alpha}}{\Gamma(2-\alpha)} - \kappa(0,z) \right) f^{\prime \prime} (t-z) \dif z.
    \nonumber
\end{eqnarray}
We deduce that the error in the approximation can be written in terms of $ \varepsilon(t) = t^{1-\alpha} / \Gamma(2-\alpha) - \kappa(0,t) $.
Hence, 
\begin{eqnarray}
    | D_t^\alpha f - \hat{D}_t^\alpha f |
    & = & 
    \left | \varepsilon(t) f^\prime (0)
    +
    \int_0^t \varepsilon(z) f^{\prime \prime} (t-z) dz \right |
    \nonumber \\
    & \le &
    | \varepsilon(t) f^\prime (0) |
    +
    \| \varepsilon \|_{0,\infty} \int_0^t  | f^{\prime \prime} (t-z) | dz
    \nonumber \\
    & \le &
    \| \epsilon \|_{0,\infty}
    \left [ 
    | f^\prime (0) | +
    \| f^{\prime \prime} \|_{0,1}
    \right ],
    \nonumber
\end{eqnarray}
where  $ \| \cdot \|_{0,\infty} $ and $ \| \cdot \|_{0,1} $ denote the $ L^\infty(0,T) $ and $ L^1(0,T) $ norms, respectively.
\end{pol}

\begin{rmk} \label{rem_l2error}
Assuming $ f^\prime(0) = 0 $ and $ f^{\prime \prime} \in L^2(0,T)$, then 
$$
    \left | D_t^\alpha f - \hat{D}_t^\alpha f \right |
    \le 
    \| \epsilon \|_{0}
    \| f^{\prime \prime} \|_{0},
$$
where $ \| \cdot \|_0 $ is the $ L^2(0,T) $ norm.
\end{rmk}

From Lemma~\ref{lem_Prony_error}, we observe that the error incurred, when re-writing the fractional derivative in terms of a finite series of Maxwell elements, is governed by $ \varepsilon $ as well as the behavior of the differentiated function $ f^{\prime\prime} $.
While for forward simulations in the time domain, the behavior of $ f $ is unknown, we can understand the impact of our approximation by directly evaluating $ \varepsilon $ over the time domain $ (0,T) $.

\subsection{Discrete Prony approximation of the Caputo derivative}
\label{meth:frac_approx_disc}
The approximation shown in Eq.~\eqref{eq_Prony_approx_series} lends itself to a straightforward discretized form by applying numerical updates to the intermediate variables $ q_k $.  
We now approximate the first derivative using backward Euler and the integral term using a midpoint approximation, i.e.
\begin{equation}\label{eq_prony_approx_final}
   \hat{\text{D}}_n^{\alpha}  f
      := 
      \frac{\beta_0}{\Delta_t} \left( f^n - f^{n-1} \right)  + \sum_{k=1}^N q_k^n, 
      \quad 
      q_k^n = e_k^2 q_k^{n-1} + e_k \beta_k \left( f^n - f^{n-1} \right),
\end{equation}
where $ e_k = \exp[ - \Delta_t / (2\tau_k)] $.
With this definition, we can examine the discretization error between the discrete operator $ \hat{\text{D}}_n^\alpha $ and its continuous counterpart $ \hat{D}_t^\alpha $, as shown in Lemma~\ref{lem_disc_Prony_error}.

\begin{lem}\label{lem_disc_Prony_error}
Let $ f: [0,T] \to \mathbb{R} $ be a real function in the time domain [0,T] and let $ f \in W^{3,\infty} (0,T) $.  
Assume the time interval $ [0,T] $ is divided into equal divisions $ \Delta_t $, giving discrete time steps $ \{ t_1, \ldots, t_{N_T} \} $ with $ t_k = k \Delta_t $ and $ T = N_T \Delta_t $.
The value $ f^n = f(t_n) $ and the value $ q_k^n $ is our approximation to $ q_k(t_n) $ based on Eq.~\eqref{eq_prony_approx_final}.
Here we assume $ q_k^0 = q_k(0) $ for all $ k = 1, \ldots, N $, and that $ \beta_k, \tau_k \in \mathbb{R}^+ $.
Then for any $ t_n \in \{ t_1, \ldots t_N \} $
\begin{equation*}
    | \hat{D}_{t_n}^\alpha f - \hat{\textnormal{D}}_n^\alpha f |  \le \Delta_t \left ( \beta_0 / 2 + C(\bs{\beta},\bs{\tau}) \Delta_t  \right )  \| f \|_{W^{3,\infty}(0,T)},
\end{equation*}
where $ C(\bs{\beta},\bs{\tau}) > 0 $ is a constant depending on the chosen $ \beta_k, \tau_k $ for all $ k = 1, \ldots, N $.
\end{lem}

\begin{pol}{\ref{lem_disc_Prony_error}}
We start by manipulating the ordinary differential equation for $ q_k $ in Eq.~\eqref{eq_Prony_approx_series}.  
Applying an integration factor and integrating over a discrete time interval $ [t_{n-1},t_n] $, we can see that 
\begin{equation*}
q_k(t_n) = \exp\left[  -\frac{\Delta_t}{\tau_k} \right] q_k(t_{n-1}) + \int_{t_{n-1}}^{t_n} \beta_k \exp\left[  \frac{s-t_n}{\tau_k} \right] f^{\prime}(s) \dif s.
\end{equation*}
Subtracting the discrete approximation shown in Eq.~\eqref{eq_prony_approx_final}, we observe that
\begin{equation} \label{eq_lem2_eq1}
  (q_k(t_n) - q_k^n) 
  = 
  \exp\left[  -\frac{\Delta_t}{\tau_k} \right] ( q_k(t_{n-1}) - q_k^{n-1} )  
  + 
  \int_{t_{n-1}}^{t_n} \beta_k \exp\left[  \frac{s-t_n}{\tau_k} \right] f^{\prime}(s) \dif s  - e_k \beta_k \left( f^n - f^{n-1} \right).
\end{equation}
We recall the result that for $ g \in W^{2,\infty}(0,T) $, one may derive the midpoint relation
\begin{equation}
  \int_{t_{n-1}}^{t_{n}} g(s) - g(t_{n-\frac{1}{2}}) \dif s 
   = 
  \int_{t_{n-1}}^{t_{n}} \int_{t_{n-\frac{1}{2}}}^{s} (s - z) g^{\prime \prime} (z) \dif z \dif s,
\end{equation}
as well as derive an inequality for the truncation error 
\begin{eqnarray}
  \int_{t_{n-1}}^{t_{n}} g(s) - g(t_{n-\frac{1}{2}}) \dif s 
  & = &
  \int_{t_{n-\frac{1}{2}}}^{t_{n}} \int_{t_{n-\frac{1}{2}}}^{s} (s - z) g^{\prime \prime} (z) \dif z \dif s
  +  
  \int_{t_{n-1}}^{t_{n-\frac{1}{2}}} \int_{s}^{t_{n-\frac{1}{2}}} (z - s) g^{\prime \prime} (z) \dif z \dif s
  \nonumber \\
  & \le &
  \| g^{\prime \prime} \|_{L^\infty(t_{n-1},t_n)} \left[ 
  \int_{t_{n-\frac{1}{2}}}^{t_{n}} \int_{t_{n-\frac{1}{2}}}^{s} | s - z | \dif z \dif s
  +  
  \int_{t_{n-1}}^{t_{n-\frac{1}{2}}} \int_{s}^{t_{n-\frac{1}{2}}} |z - s|  \dif z \dif s
  \right]
  \nonumber \\
  & \le &
  \frac{\Delta_t^3}{24} \| g^{\prime \prime} \|_{L^\infty(t_{n-1},t_n)}.
  \nonumber 
\end{eqnarray}
Utilizing the midpoint relation, the second half of Eq.~\eqref{eq_lem2_eq1} can be written as
\begin{eqnarray}
   &  &
   \int_{t_{n-1}}^{t_n} \beta_k \exp\left[  \frac{s-t_n}{\tau_k} \right] f^{\prime}(s) \dif s  - e_k \beta_k \left( f^n - f^{n-1} \right) 
       \nonumber \\
       & & \hspace{10mm} 
       =
        \int_{t_{n-1}}^{t_n} \beta_k \exp\left[  \frac{s-t_n}{\tau_k} \right] f^{\prime}(s)  - \beta_k e_k f^{\prime}(t_{n-\frac{1}{2}}) \dif s  
       - e_k \beta_k \left( f^n - f^{n-1}  -  \Delta_t  f^{\prime}(t_{n-\frac{1}{2}})  \right) 
       \nonumber \\
       & & \hspace{10mm} 
       =
        \int_{t_{n-1}}^{t_n} \beta_k \exp\left[  \frac{s-t_n}{\tau_k} \right] f^{\prime}(s)  - \beta_k e_k f^{\prime}(t_{n-\frac{1}{2}}) \dif s  
       - e_k \beta_k \int_{t_{n-1}}^{t_n} f^\prime(s) - f^{\prime}(t_{n-\frac{1}{2}}) \dif s.  
       \nonumber \\
\end{eqnarray}
Note that each integral represents the truncation error of a midpoint approximation, i.e.
\begin{eqnarray}
   &  &
   \left | \int_{t_{n-1}}^{t_n} \beta_k \exp\left[  \frac{s-t_n}{\tau_k} \right] f^{\prime}(s) \dif s  - e_k \beta_k \left( f^n - f^{n-1} \right) \right|
       \nonumber \\
       & & \hspace{10mm} 
       \le
       \frac{\beta_k \Delta_t^3}{24} \left \| \left ( \exp \left [ \frac{s-t_n}{\tau_k} \right ] f^\prime (s) \right)^{\prime \prime} \right \|_{L^\infty(t_{n-1},t_n)} 
       + \frac{e_k \beta_k \Delta_t^3}{24} \left \| f^{\prime \prime \prime} \right \|_{L^\infty(t_{n-1},t_n)} 
       \nonumber \\
       & & \hspace{10mm} 
       \le
        \frac{\beta_k \Delta_t^3}{24} \left ( \frac{1}{\tau_k^2} \| f^\prime \|_{L^\infty(t_{n-1},t_n)} + \frac{2}{\tau_k} \| f^{\prime \prime} \|_{L^\infty(t_{n-1},t_n)} + (1+e_k) \| f^{\prime \prime \prime} \|_{L^\infty(t_{n-1},t_n)} \right )
       \nonumber \\
       \nonumber \\
       & & \hspace{10mm} 
       \le
        \Delta_t^3 C_0(\beta_k,\tau_k) \| f \|_{W^{3,\infty}(t_{n-1},t_n)},
       \nonumber
\end{eqnarray}
where $ C_0(\beta_k,\tau_k) = (\beta_k / 24) \max \{ \tau_k^{-2}, \tau_k^{-1}, (1+e_k) \} $.  
Inserting the truncation error into Eq.~\eqref{eq_lem2_eq1}, we can see that
\begin{eqnarray}
  | q_k(t_n) - q_k^n |
  & \le &
  e_k^2 | q_k(t_{n-1}) - q_k^{n-1} | 
  + 
  \left |  \int_{t_{n-1}}^{t_n} \beta_k \exp\left[  \frac{s-t_n}{\tau_k} \right] f^{\prime}(s) ds  - e_k \beta_k \left( f^n - f^{n-1} \right) \right |
  \nonumber \\
  & \le &
  e_k^2 | q_k(t_{n-1}) - q_k^{n-1} | 
  + 
  \Delta_t^3 C_0(\beta_k,\tau_k) \| f \|_{W^{3,\infty}(t_{n-1},t_n)}
  \nonumber \\
  & \le &
  e_k^2 | q_k(0) - q_k^{0} | 
  + 
  C_0(\beta_k,\tau_k) \Delta_t^3  \sum_{k=1}^{n} \| f \|_{W^{3,\infty}(t_{k-1},t_k)}
  \nonumber \\
  & \le &
  C_0(\beta_k,\tau_k) T \Delta_t^2  \| f \|_{W^{3,\infty}(0,T)},
  \nonumber 
\end{eqnarray}
making the update formula for $ q_k $  $ \mathcal{O} (\Delta_t^2) $.  
Subtracting the definition of $ \hat{D}_t^\alpha $ from $ \hat{\text{D}}_n^\alpha $, we can see that 
\begin{eqnarray}
   | \hat{D}_{t_n}^\alpha f - \hat{\text{D}}_n^\alpha f |  
   & = &
   \left |
   \beta_0 \left ( f^\prime(t_n) - \Delta_t^{-1} ( f^n - f^{n-1}) \right ) + \sum_{k=1}^N \left( q_k(t_n) - q_k^n \right) 
   \right |
   \nonumber \\
   & \le &
   \beta_0 | f^\prime(t_n) - \Delta_t^{-1} ( f^n - f^{n-1}) | + \sum_{k=1}^N | q_k(t_n) - q_k^n |
   \nonumber \\
   & \le &
   \beta_0 \left | \frac{1}{\Delta_t} \int_{t_{n-1}}^{t_n} (s-t_{n-1}) f^{\prime \prime}(s) ds  \right | 
   + 
   C(\bs{\beta},\bs{\tau}) \Delta_t^2 \| f \|_{W^{3,\infty}(0,T)}
   \nonumber \\
   & \le &
   \frac{\beta_0 \Delta_t}{2}   \| f^{\prime \prime} \|_{L^{\infty}(t_{n-1},t_n)}
   + 
   C(\bs{\beta},\bs{\tau}) \Delta_t^2 \| f \|_{W^{3,\infty}(0,T)},
\end{eqnarray}
where $ C(\bs{\beta},\bs{\tau}) = \sum_{k=1}^N C_0(\beta_k,\tau_k)  $.  Generalizing for any time point completes the proof.
\end{pol}

From Lemma~\ref{lem_disc_Prony_error} we observe the typical discretization errors expected for backward Euler and midpoint based error estimates.  
The result shows that the discrete approximation converges to the approximate Prony series approximation with $ \mathcal{O} (\Delta_t ) $.
Using backward Euler for the first-order derivative term is not necessary, and higher-order approximations could easily be used in Eq.~\eqref{eq_prony_approx_final}.
Here we select this form for ease and, in the context of biomechanics simulations, for stability considerations.
Combining Lemmas~\ref{lem_Prony_error} and~\ref{lem_disc_Prony_error}, we can derive an error estimate for this Prony series-based approximation in Theorem~\ref{thm_prony_error}.

\begin{thm} \label{thm_prony_error}
\sloppy Given the assumptions of Lemma~\ref{lem_disc_Prony_error}, the error between $ \hat{\textnormal{D}}_n^\alpha f $ and $ D_{t_n}^\alpha f $ for any $ t_n \in \{ t_1, \ldots, t_{N_T} \} $ is given by
$$
    | \hat{\textnormal{D}}_n^\alpha f - D_{t_n}^\alpha f | \le 
     \| \epsilon \|_{0,\infty}
    \left [ 
    | f^\prime (0) | +
    \| f^{\prime \prime} \|_{0,1}
    \right ]
    +
     \Delta_t \left ( \beta_0 / 2 + C(\bs{\beta},\bs{\tau}) \Delta_t  \right )  \| f \|_{W^{3,\infty}(0,T)}.
$$
\end{thm}
\begin{pot}{\ref{thm_prony_error}}
This is a straightforward result of the triangle inequality, Lemma~\ref{lem_Prony_error} and~\ref{lem_disc_Prony_error}.
\end{pot}

Based on this error estimate, we can see that a refinement with $ \Delta_t $ eliminates the portion of error responsible for discretization, while the appropriate selection of the positive constants used in Eq.~\eqref{eq_prony_approx_final} can reduce the Prony approximation error given by the truncation $ \varepsilon $ (as we will later show).  
We note that decreasing $ \Delta_t $ will eventually lead to a plateau in the error convergence based on how well the Prony series approximates the true fractional derivative.
However, increasing the number of Prony terms $N$ elicits a more complex response.
For large $ \Delta_t $, the error can become dependent on the second term in Theorem~\ref{thm_prony_error} and, thus, the constant $ C (\bs{\beta},\bs{\tau}) $.
This constant, in general, grows with more Prony terms and can have truncation scale factors that grow larger with larger $ N $.
However, for sufficiently fine $ \Delta_t $, increasing the number of terms shifts the error onto $ \varepsilon $ and thus improves like the truncation error. 

\subsection{Optimizing the Numerical approximation $\hat{D}_t^{\alpha} $} \label{meth:frac_opt}

The approximation introduced in Eq.~\eqref{eq_Prony_approx} relies on the effective identification of the parameters $ \bs{\theta}_\alpha $ that should be chosen to ensure optimal accuracy.
Ideally, $ \bs{\theta}_\alpha $ would be selected so as to minimize the error between the true Caputo derivative and the approximation in some suitable norm.
For example, if we consider the $L^2(0,T) $ norm, then we would seek to minimize the squared error 
\begin{equation}
F(f,\bs{\theta}_\alpha) = \frac{1}{2} \int_0^T  [ \hat{D}_t^\alpha f - D_t^\alpha f ]^2 \; \dif t,
\end{equation}
where $ F (f, \bs{\theta}_\alpha) $ gives the squared $ L^2 $-norm error between the approximate and the true Caputo derivatives for the given function $ f $ and set of parameters $ \bs{\theta}_\alpha $.  
While the parameters $ \bs{\theta}_\alpha $ could be chosen to minimize $ F(f, \bs{\theta}_\alpha ) $ this approach presents challenges as we often do not know $ f $ prior to simulation.
Hence, instead, we must select our approximation to be suitably valid for all $ f $, \emph{i.e.} we find $ \bs{\theta}_\alpha $ such that (for $ 0 < \alpha < 1 $),
\begin{equation} \label{eq_min_principle}
F(f,\bs{\theta}_\alpha,\alpha) := \left \{ \min F(f,\bs{\varphi},\alpha), \; \bs{\varphi} \in \mathbb{R}^{2N+1} \right \}, \quad  \text{ for any } f \in L^\infty(0,T) .
\end{equation}
The problem posed in Eq.~\eqref{eq_min_principle} is generally challenging to solve due to the flexibility in the behavior of $ f $.
In order to approximately minimize $ F$, suppose we instead consider only those functions $ f \in \mathcal{F}(M) $ where
\begin{equation}
\mathcal{F}(M) = \left \{ f \in L^2(0,T)  \; \Big | \;  f(t) = \text{Re} \left \{ \dsum_{k=1}^M C_k \exp ( i \omega_k t ) \right \}, \; \text{for some} \; C_0, \ldots, C_M  \in \mathbb{E} \right \}
\end{equation}
represents the set of functions that can be written using a Fourier basis~\cite{tseng2000computation}.  Here we note that $ \omega_k = k \omega^\star $, where $ \omega^\star = 2 \pi / T $ is the base frequency of $ f $ on the time interval $ [0,T] $.  In this context, our function $f $ and its fractional derivative (for $ t >> 0 $) can be written as
\begin{equation}
f(t) = \text{Re} \left \{ \sum_{k=1}^M C_k \exp ( i \omega_k t ) \right \},
\quad
D_t^\alpha f(t) \approx  \text{Re} \left \{ \sum_{k=1}^M C_k (i \omega_k)^\alpha  \exp ( i \omega_k t ) \right \}.
\end{equation}
In contrast, the approximate fractional derivative $ \hat{D}_t^\alpha $ based on Eq.~\eqref{eq_Prony_approx} can be written as,
\begin{equation}
\hat{D}_t^\alpha f(t) =  \text{Re} \left \{  \sum_{k=1}^M C_k \left( \beta_0 i \omega_k + \sum_{m=1}^N \beta_m (\tau_m \omega_k) \frac{(\tau_m \omega_k + i)}{(\tau_m \omega_k)^2+1}  \right) \exp \{ i \omega_k t \} \right \}.
\end{equation}
The difference between the true and approximate operators is then simply,
\begin{equation}
\hat{D}_t^\alpha f(t)  - D_t^\alpha f(t) \approx 
                \text{Re} \left \{ \sum_{k=1}^M C_k \exp ( i \omega_k t )  \left [ \beta_0 i \omega_k + \sum_{m=1}^N \beta_m (\tau_m \omega_k) \frac{(\tau_m \omega_k + i)}{(\tau_m \omega_k)^2+1} - (i \omega_k)^\alpha  \right] \right \},
\end{equation}
leaving the operator error dependent on the size of the complex number given in the brackets.  
This means that for us to obtain a good approximation for any $ f \in \mathcal{F}(M) $, we must choose $ \bs{\theta}_\alpha $ such that we minimize the bracketed term.
Hence, we may reduce the minimization problem in Eq.~\eqref{eq_min_principle} to an analogous normalised system where
\begin{equation} \label{eq_algebraic_constraints}
\frac{1}{k^\alpha} \sum_{m=1}^N \hat{\beta}_m \frac{(k \hat{\tau}_m)^2}{(k \hat{\tau}_m)^2+1} = \cos \left( \frac{\pi \alpha}{2} \right),
\quad
\hat{\beta}_0 k^{1-\alpha} 
-  \frac{1}{k^\alpha} \sum_{m=1}^N \hat{\beta}_m \frac{k \hat{\tau}_m}{(k \hat{\tau}_m)^2+1} = \sin \left( \frac{\pi \alpha}{2} \right), 
\end{equation}
for $ k = 1, \ldots M $.  
Here, we note that parameters were normalized (with hats) to enable independence from the base frequency $ \omega^\star $, with $  \beta_0 =  \hat{\beta}_0 (\omega^\star)^{\alpha-1} $, $  \beta_m = \hat{\beta}_m  (\omega^\star)^\alpha $ and $  \tau_m = \hat{\tau}_m  / \omega^\star $ ($ m = 1, \ldots N $). 
This leads to $2M $ constraint equations that, if hold, imply equivalence of the true and approximate Caputo derivative.  
Assuming that our function is composed of many more Fourier terms than the number of Prony terms (i.e. $ M >> N $), Eq.~\eqref{eq_algebraic_constraints} leads to a nonlinear over-constrained system that can be solved by optimization.

\subsection{Numerical implementation} \label{meth:implementation}
    Resulting parameters for approximating different fractional derivatives can be found in the MATLAB~\cite{matlab2018version} codes provided as supplementary material.
    Here optimized parameters can be returned for values of $ N \in [3,15] $ and $ \alpha \in [0,1] $ for a specified base frequency $ \omega^\star $.
    Codes and examples are also provided showing how the approximation can be utilized for approximating the Caputo derivative.
    Other methods, such as the midpoint (MP) or Gr{\"{u}}nwald-Letnikov (GL) approximations, are also provided.


\section{Large Deformation Viscoelastic Solid Mechanics} \label{meth:solid_mech}
%

\subsection{Review of kinematics and kinetics}\label{sect:kinematics}

In this section, we briefly review the kinematics and kinetics for hyperelastic materials.
Let $\Omega_0 \subset \mathbb{R}^3 $ denote the reference configuration of a three-dimensional solid (coordinates $ \bs{X} $) and $ \Gamma_0 $ denote its boundary.
The boundary is split into groups $ \Gamma_0 = \Gamma_0^N \cup \Gamma_0^D $ based on either loading (subscript $ N $) or displacement boundary conditions are applied (subscript $ D $).
Over a time interval $t\in [0,T] $, its' configuration is $ \Omega(t) \subset \mathbb{R}^3 $, with physical coordinates $\bs{x}$.
The coordinates are related by
\begin{equation}
    \bs{x}(t) = \bs{u}(t) + \bs{X}, 
    \qquad
    x_i (t) = u_i(t) + X_i.
    \label{eq_x_to_X}
\end{equation}
This mapping in Eq.~\eqref{eq_x_to_X} is assumed to be bijective with $ \bs{u} : \Omega_0 \times [0,T] \to \mathbb{R}^3 $ denoting the material displacements.  
Note that $\bs{X}$ and $\bs{x}$ are the vector positions, with components $X_i$ and $x_i$. 
The local deformation gradient $\bs{F}$ defines the local mapping of infinitesimal vectors, which is given by 
\begin{equation*}
    \bs{F} = \dpd{\bs{x}}{\bs{X}} = \dpd{\bs{u}}{\bs{X}} + \bs{I}, 
    \qquad 
    F_{ij} = \frac{ \partial u_i }{ \partial X_j} + \delta_{ij},
\end{equation*}
where $\bs{I}$ is the identity tensor.
The determinant $ J = \det \bs{F} $ defines the relative change in volumes, providing key information about material growth / decay.
The right and left Cauchy-Green tensors~\cite{lai2009introduction}, defining the material stretch, are
\begin{equation*}
    \bs{C} = \bs{F}^{T} \bs{F}, 
    \qquad 
    \bs{B}=\bs{F}\bs{F}^{T}. 
\end{equation*}
The deformation gradient and the left and right Cauchy-Green tensors comprise a common set of measures used to characterize the constitutive behavior of a material.  
For a hyperelastic material, the strain-energy density function is given by $ \Psi = \Psi(\bs{C})$, which defines the stored energy resulting from material deformations.  
From the second law of thermodynamics, the strain-energy function can be related to the second Piola-Kirchhoff stress tensor (PK2) as
\begin{equation}
    \bs{S} = 2 \dpd{\Psi(\bs{C})}{\bs{C}}, \qquad S_{ij} = 2 \frac{\partial \Psi(C_{ij})}{\partial C_{ij}},
    \label{eq_pk2}
\end{equation}
which, in turn, is related to the first Piola-Kirchhoff and Cauchy stress tensors ($\bs{P}$ and $ \bs{\sigma} $, respectively) by push forward operations on PK2.  In this case,
\begin{equation}
    \bs{\sigma} = \frac{1}{J} \bs{F} \bs{S} \bs{F}^T = \frac{1}{J} \bs{P} \bs{F}^T.
    \label{eq_cauchystress}
\end{equation}
For isotropic biomechanical materials, the strain-energy function is often defined using invariants of $ \bs{C} $~\cite{Bonet1997}, 
\begin{equation}
    \I = \bs{C} : \bs{I}, \quad \II = \bs{C}:\bs{C}, \quad \III = \det \bs{C} 
\end{equation}
with $ \Psi = \Psi (\I,\II, \III) $.  
For incompressible materials, the third invariant $ \III $ is constrained to $ 1 $ by splitting the strain-energy function into a deviatoric ($\Psi_d$) and a volumetric components with $ \Psi = \Psi_d + p (\III^{1/2}-1) $, where $ p : \Omega \times [0,T] \to \mathbb{R} $ denotes the hydrostatic pressure.
Often the strain-energy function takes exponential forms, as seen in the Fung model \cite{fung1993biomechanics}.  
However, with many tissues exhibiting anisotropic material response \cite{holzapfel2002nonlinear}, additional invariants have been introduced to denote stretch with respect to specific microstructural directions.  
These laws can be naturally extended, but for the purposes of this work, we will focus on more simple isotropic forms.

\subsection{Extension of hyperelastic to viscoelastic models}\label{meth:frac_model}

For constitutive modeling of viscoelastic behavior, a number of model approaches have been presented in the literature.
Fig.~\ref{fig:maxandkelvinmodels} illustrates some of the common models introduced, particularly for linear viscoelastic materials, using spring, dashpot or spring-pot elements.
Extension of these forms toward nonlinear material models has become more common, with models following the analogy presented in Fig.~\ref{fig:maxandkelvinmodels}.
For example, a nonlinear viscoelastic Maxwell model (Fig.~\ref{fig:maxandkelvinmodels}A) could be written as 
\begin{equation}
    \bs{S} + \frac{\eta_1}{E_1}\dot{\bs{S}} = \eta_1 \dot{\bs{S}}_v
\end{equation}
and the nonlinear viscoelastic Kelvin-Voigt model (Fig.~\ref{fig:maxandkelvinmodels}B) as,
\begin{equation}
    \bs{S} = E_0 \bs{S}_e + \eta_1 \dot{\bs{S}}_v.
\end{equation}
Here $ \bs{S}_e $ denotes a symmetric tensor denoting the hyperelastic response of the material, while $ \bs{S}_v $ denotes the viscoelastic response. 
Note that these stresses denote the non-hydrostatic components, which can be included based on the application.
Commonly, both take a similar form to hyperelastic materials commonly encountered in the literature.
These models provide an exponential decay in the viscoelastic phenomena, which can yield viscoelastic response times that are incompatible with experimental data.
Other formulations, such as introduced by Schapery et al. \cite{schapery1966theory}, provide a more generalized view for characterizing material response.
Another example is the model developed by Freed and Diethelm \cite{freed2006fractional}, i.e.
\begin{equation}
    \bs{S}(t) = \left(\mu_\infty +(\mu_0-\mu_\infty)\right)G(t)\dpd{\Psi(\bs{E})}{\bs{E}} + (\mu_0-\mu_\infty) \int_0^t M(t-s)\dpd{\Psi(\bs{E})}{\bs{E}}\dif s,
\end{equation}
    where $G(t)$ is a relaxation function and $M(t)$ is a memory function for the history dependency. 
    
    An alternative approach, that has been successful in the experimental literature, is use of fractional elements~\cite{magin2006fractional,holm2019waves}.
    This is due to the progressive memory that these models exhibit that follows a power spectrum commonly observed in biological tissues for moderate frequencies.
    In this case, we can envision development of nonlinear fractional viscoelastic models, such as the fractional Kelvin Voigt model (Fig.~\ref{fig:maxandkelvinmodels}E), \emph{e.g.}
\begin{equation}
    \bs{S} = \bs{S}_e + D_t^\alpha \bs{S}_v .
\end{equation}
    Here $\alpha = 0$ reduces to a pure elastic response and $\alpha = 1$ converts the latter part of the equation to a purely viscous-like response. 
    The advantage of this approach is the generality of its application and ability to capture the dynamic viscoelastic behavior of tissues with the addition of the fractional differential operator.
    
\subsection{Conservation of mass and momentum}

With an appropriately defined constitutive model, the conservation of mass and momentum for an incompressible solid material is given by~\cite{holzapfel2002nonlinear},
\begin{eqnarray}
\partial_t (\varrho J \bs{v}) - \nabla_{\bs{X}} \cdot (\bs{F} \bs{S}) & = & J \bs{b}, 
        \quad \text{ in } \Omega_0 \times [0,T],
        \label{eq_conservation_momentum} \\ 
J - 1 & = & 0, \quad \text{ in } \Omega_0 \times [0,T], 
        \label{eq_conservation_mass}
\end{eqnarray}
where $ \partial_t $ is the Lagrangian derivative (taken with respect to fixed coordinates $ \bs{X} $),
$ \varrho $ is the material density, $ \bs{v} = \partial_t \bs{u} $ is the material velocity, and $ [0,T] $ is the time interval of interest.
The mechanical system is then subject to an initial condition and boundary conditions, e.g.
\begin{eqnarray}
\bs{u}(\cdot,0) =  \bs{u}_0, 
& \quad &
\bs{v}(\cdot,0)  =  \bs{v}_0 
        \quad \text{ in } \Omega_0,
        \label{eq_dispvel_ic} \\ 
\bs{v}(\cdot,t)  =  \bs{v}_D,
        \quad \text{ on } \Gamma_0^D \times [0,T], 
& \quad &        
(\bs{F} \bs{S}) \cdot \bs{N} = \bs{t}_N 
        \quad \text{ on } \Gamma_0^N \times [0,T],
        \label{eq_veltraction_bc}
\end{eqnarray}
where $ \bs{u}_0  , \bs{v}_0 : \Omega_0 \to \mathbb{R}^3  $ are given initial conditions, $ \bs{v}_D : \Gamma_0^D \to \mathbb{R}^3 $ given information about the rate of deformation on the boundary $ \Gamma_0^D  $, and $ \bs{t}_N : \Gamma_0^N \to \mathbb{R}^3 $ is given information about the applied traction on boundary $ \Gamma_0^N $.
All examples considered herein assume that no history of deformation is present in the material for $ t < 0 $ (e.g., the fractional derivative is defined w.r.t. $ t = 0$).

\section{Finite Element Approximation for Viscoelastic Solid Models} \label{meth:solid_mech_num}

This section reviews the integration of the nonlinear fractional viscoelastic model into a finite element solid mechanics framework.
For development and discretization of the weakform system, we follow the typical mixed finite element $ \bs{u}-p $ formulation~\cite{zienkiewicz2000finite, hughes2012finite} (section~\ref{meth:solid_mech_wf}) and extend it to account for the discrete fractional order model (section~\ref{meth:solid_mech_stress_disc}).
A summary of the implementation steps is provided in section~\ref{meth:solid_mech_solproc}.
The linearized system is introduced in section~\ref{meth:solid_mech_discwf_stability}, and stability estimates are proven for the discrete problem.

\subsection{Continuous and Discrete Weakform System} \label{meth:solid_mech_wf} 

The weak formulation for the mechanical system follows the typical procedure.
Multiplying Eqs.~\eqref{eq_conservation_momentum}-\eqref{eq_conservation_mass} by test functions $ \bs{w} $ and $ q $, respectively, integrating and applying integration by parts, the weak problem can be written as: find $ (\bs{u} (t), \bs{v}(t), p(t)) \in \bs{\mathcal{U}} \times \bs{\mathcal{V}}_D \times \mathcal{P} $ such that,
\begin{equation} \label{eq_fem_wf}
    \left( \partial_t ( \varrho \bs{v} J) , \bs{y}  \right)
    +
    (\bs{F} \bs{S}, \nabla_{\bs{X}} \bs{y} )
    - 
    \big < \bs{t}_N, \bs{y}  \big >
    -
    (J\bs{b}, \bs{y})
    +
    (q,J-1)
    + 
    (\bs{v} - \partial_t \bs{u}, \bs{w} )
    = 
    0, 
\end{equation}
for all $ (\bs{w}, \bs{y}, q) \in \bs{\mathcal{U}} \times \bs{\mathcal{V}}_0 \times \mathcal{P} $ and all times $ t \in [0,T] $. Here, $ \bs{\mathcal{U}} $, $ \bs{\mathcal{V}} $ and $ \mathcal{P} $ denote appropriate spaces~\cite{asner2017patient}, their subscripts $ D $ and $ 0 $ denote spaces equipped with appropriate Dirichlet or zero Dirichlet conditions on $ \Gamma_0^D $, and the operators $ ( \cdot, \cdot) $ and $ \big < \cdot ,  \cdot  \big > $ denote the inner product on $ \Omega_0 $ and $ \Gamma_0^N $, respectively.

A discrete-in-time discretization of Eq.~\eqref{eq_fem_wf} can be considered using Backward Euler for temporal derivative approximations.
Here we define a series of time steps $ \{ t_n \}_{n=0}^{N_T} $ (as seen in section~\ref{sect:current_frac_meths}), and letting superscript $ n $ denote a variable quantity at time $ t_n $, and define the Backward temporal operator, $ \delta_t f^n = (f^n - f^{n-1}) / \Delta_t$.
Defining appropriate inf-sup compatible discrete spaces, denoted with the subscript $ h $ (in later results, we use $ \mathbb{Q}^2 - \mathbb{Q}^2 - \mathbb{Q}^1$ Taylor-Hood elements~\cite{brezzi1991stability}), the discrete weak form of Eq.~\eqref{eq_fem_wf} can be written as: find $ (\bs{u}_h^n, \bs{v}_h^n, p_h^n ) \in \bs{\mathcal{U}}^h \times \bs{\mathcal{V}}_D^h \times \mathcal{P}^h $ such that,
\begin{equation} \label{eq_fem_disc_wf}
    \left ( \delta_t (\varrho \bs{v}_h^n J^n), \bs{y}_h \right )
    +
    (\bs{F}^n \bs{S}^n, \nabla_{\bs{X}} \bs{y}_h )
    - 
    \big < \bs{t}_N, \bs{y}_h  \big >
    +
    (q_h,J^n-1)
    -
    (J^n\bs{b}^n, \bs{y}_h)
    +
    (\bs{v}_h^n - \delta_t \bs{u}_h^n, \bs{w}_h )
    = 
    0, 
\end{equation}
for all $ (\bs{w}_h, \bs{y}_h, q_h) \in \bs{\mathcal{U}}^h \times \bs{\mathcal{V}}_0^h \times \mathcal{P}^h $ and all $ n = 1, \ldots N_T $.  
We note that, under the scheme presented, the solution for $ \bs{u}_h^n $ can be implicitly defined based on $ \bs{v}_h^n $ as
$$ \bs{v}_h^n = \delta_t \bs{u}_h^n $$
at node points.

\subsection{Stress Discretization and Integration with the Prony Approximation} \label{meth:solid_mech_stress_disc}

Evaluation of the discrete finite element equation in Eq.~\eqref{eq_fem_disc_wf} follows the typical procedure applied in mechanics.
Namely, the weighted residual equations are evaluated using an appropriate form of quadrature based on the element type and computational complexity of the integrand.  
In the context of fractional viscoelastic models, as introduced in section~\ref{meth:frac_model}, the only additional challenge is determining a scheme for handling the viscoelastic stress tensor.  
Using the prony approximation introduced in Eq.~\ref{eq_prony_approx_final}, the fractional Kelvin-Voigt model can be written as,
\begin{equation} \label{eq_disc_frac_stress}
    \bs{S}^n = \bs{S}_e^n 
             + \beta_0 \delta_t \bs{S}_v^n
             + \sum_{k=1}^N \bs{Q}_k^n
             + p_h^n \bs{S}_p^n,
\end{equation}
where $ \bs{S}_e^n = \bs{S}_e ( \bs{C}^n) $, $ \bs{S}_v^n = \bs{S}_v (\bs{C}^n) $, $ \bs{S}_p^n = J^n [ \bs{C}^n ]^{-1} $ and $ \bs{Q}_k^n $ is defined by the update formula,
\begin{equation*}
    \bs{Q}_k^n = e_k^2 \bs{Q}_k^{n-1} + \beta_k e_k \Delta_t \delta_t \bs{S}_v^n.
\end{equation*}
The approximation for the current value of the second Piola Kirchhoff stress given in Eq.~\eqref{eq_disc_frac_stress} introduces the intermediate variables, $ \bs{Q}_k $, representing the evolution of the material history due to different Prony terms.  
Within a finite element discretization, Eq.~\eqref{eq_disc_frac_stress} could be realized weakly, interpolating $ \bs{Q}_k $ using basis functions.
However, this form does not ensure equivalence to the original form observed in Eq.~\eqref{eq_fem_disc_wf}.
Instead, the values of $ \bs{Q}_k $ are \emph{discretized} at quadrature points where the residual equation (and, consequently, the stress) must be evaluated.
Fig.~\ref{fig:fe_discretization} illustrates this strategy for $ \mathbb{P}^2-\mathbb{P}^1 $ mixed-problem set on triangles using an appropriate quadrature rule~\cite{lyness1975moderate}.  In Fig.~\ref{fig:fe_discretization}B we see the mixed element showing node points for displacement and pressure and in Fig.~\ref{fig:fe_discretization}C the quadrature points at which values of $ \bs{Q}_k $ are stored.

For practical implementation purposes, it is convenient to re-arrange Eq.~\eqref{eq_disc_frac_stress} to consolidate terms based on their dependence on the current time step and past time step, \emph{e.g.}
\begin{equation} \label{eq_disc_frac_stress2}
   \bs{S}^n 
      = 
      \bs{S}_e^n + \gamma \bs{S}_v^n  + p_h^n \bs{S}_p^n
      - \gamma \bs{S}_v^{n-1} + \sum_{k=1}^N e_k^2 \bs{Q}_k^{n-1},
\end{equation}
where $ \gamma = \beta_0 / \Delta_t + \sum_{k=1}^N \beta_k e_k $.  
From Eq.~\eqref{eq_disc_frac_stress2}, the first three terms depend on $ n $ while the later two terms depend on variables at $ n-1 $.
As a result, to iteratively solve Eq.~\eqref{eq_fem_disc_wf}, one needs to evaluate the last two terms once per time step.
While the first three terms vary each iteration, the complexity of computing these terms over a standard hyperelastic formulation is the addition of $ \gamma \bs{S}_v^n$.

\subsection{Solution procedure} \label{meth:solid_mech_solproc}

The computational model framework presented was implemented into CHeart \cite{lee2016multiphysics}, a custom multiphysics finite element solver. 
The solution procedure follows the Shamanskii-Newton-Raphson (SNR) method~\cite{shamanskii1967modification} shown in Algorithm~\ref{monolithic-partitioned-alg}.
Due to the formulation devised in section~\ref{meth:solid_mech_wf}, the final linearized system is written for variables $ \bs{v}_h^n $ and $ p_h^n $ with the current displacement calculated as discussed previously.
In this case, our residual functional can be written as,
\begin{equation} \label{eq_res_eqn}
    R(\bs{v}_h^n, p_h^n; \bs{w}_h, q_h) :=  \left (\varrho \delta_t (\bs{v}_h^n J^n ), \bs{w}_h \right )
    +
    (\bs{F}^n \bs{S}^n, \nabla_{\bs{X}} \bs{w}_h )
    - 
    \big < \bs{t}_N, \bs{w}_h  \big >
    +
    (q_h,J^n-1)
    -
    (J^n\bs{b}^n, \bs{w}_h).
\end{equation}
Owing to the fact that $ \bs{v}_h^n $ can be expressed as the weighted sum of vector basis functions (and similar for $ p_h^n $, but using scalar basis functions),
Eq.~\eqref{eq_res_eqn} can be written as 
$$ R(X^n; Y) = R(\varPhi_{\bs{v}} V^n, \varPhi_p P^n; \varPhi_{\bs{v}} W, \varPhi_p Q ) $$
where $ X^n = (V^n, P^n) $ denotes the coefficients of our solution and $ Y = (W, Q) $ the weights of our test functions. 
Finally, the vector residual function can be expressed as the gradient with respect to $ Y $, \emph{e.g.}
\begin{equation}
  \bs{R}(U^n) = \nabla_Y R(U^n, Y) = \bs{0}.
\end{equation}

\begin{algorithm}
    \caption{~~\texttt{Viscoelastic Mechanics Algorithm.}}\label{monolithic-partitioned-alg}
    \begin{algorithmic}[1]
        \State Given initial condition $X^{0}$ (and $ \bs{Q}^0$ at all quadrature points).
        \State \texttt{Compute} $~~\boldsymbol{J}_{\beta} = \nabla_X~\boldsymbol{R} (X^0)$, $~[\boldsymbol{J}_{\beta}]^{-1}$.
        \State
        \For{~($n = 1:N$)~}
            \State
            \State \texttt{Set} $~~X^{n, 0} = X^{n-1}$.
            \State \texttt{Compute} $~~\bs{R} (X^{n,0})$.
            \State \texttt{Compute} $~~r = \| \bs{R}(X^{n,0} ) \|$.
            \State \texttt{Set} $~~k = 0$.
            \State
            \While{~($r > TOL$)~}
                \State
                \State \texttt{Compute} $~~\delta X~~= -[\boldsymbol{J}_{\beta}]^{-1}~\boldsymbol{R}~(X^{n,k})$. 
                \State \texttt{Find} $~~\min_{\alpha \in [0, 1]} \| \boldsymbol{R} (X^{n,k} + \alpha~\delta X ) \|$.
                \State \texttt{Set} $~~X^{n,k+1}~~= X^{n,k} +  \alpha~\delta X$.
                \If{~($\| \boldsymbol{R} (X^{n,k+1}) \| > \gamma~r$~~or~~$k > ITER$)~}
                    \State \texttt{Compute} $~~\boldsymbol{J}_{\beta} = \nabla_{X}~\boldsymbol{R}~(X^{n,k+1})$, $~[\boldsymbol{J}_{\beta}]^{-1}$.
                \EndIf
                \State \texttt{Compute} $~~r = \| \boldsymbol{R}~(X^{n,k+1}) \|$.
                \State \texttt{Set} $~~k = k + 1$.
            \EndWhile
            \State \texttt{Update} $ \bs{Q}^n = e_k^2 \bs{Q}_k^{n-1} + \beta_k e_k \Delta_t \delta_t \bs{S}_v^n $ at quadrature points.
        \EndFor
    \end{algorithmic}
\end{algorithm}

\subsection{Stability of the Discrete Weak form System}\label{meth:solid_mech_discwf_stability}

In this section, stability is analyzed for linear incompressible fractional viscoelastic materials shown in Eq.~\eqref{eq_conservation_momentum_lin} and~\eqref{eq_conservation_momentum_lin}.
The problem depends on parameters $ \varrho, E, \eta \in \mathbb{R}^+ $ which, for ease, we take as positive spatiotemporal constants denoting the material density, elastic stiffness and fractional viscoelastic stiffness, respectively.
Here $ \bs{b} $ denotes body forces and $ \bs{u}_0  $ and $ \bs{v}_0  $ are given initial conditions.
All examples considered herein assume that no history of deformation (\emph{e.g.} the fractional derivative is defined w.r.t. $ t = 0$).
\begin{eqnarray}
\varrho \partial_t (\bs{v}) - \nabla_{\bs{X}} \cdot \bs{\Sigma} = \bs{b}, \; \;
        & \quad & \text{ in } \Omega_0 
        \label{eq_conservation_momentum_lin} \\ 
\nabla_{\bs{X}} \cdot \bs{v} = 0, \; \;
        & \quad & \text{ in } \Omega_0 
        \label{eq_conservation_mass_lin} \\
\bs{u}(\cdot,0)  =  \bs{u}_0, \qquad
\bs{v}(\cdot,0)  =  \bs{v}_0, 
        & \quad & \text{ in } \Omega_0 
        \label{eq_vel_ic_lin} \\ 
\bs{v} = \bs{0}  
        \quad \quad
        \text{ on } \Gamma_0^D, 
        \quad 
\bs{\Sigma} \cdot \bs{N} = \bs{0}, \; \;
        & \quad & \text{ on } \Gamma_0^N 
        \label{eq_bc_lin} 
\end{eqnarray}
The constitutive model for the linear elastic material is given by
\begin{equation} \label{eq_lin_constitutive}
    \bs{\Sigma} = E \bs{D} \bs{u} + \eta D_t^\alpha (\bs{D} \bs{u}) + p \bs{I}.
\end{equation}
with $ \bs{D}(\cdot) = \frac{1}{2} \left ( \nabla_{\bs{X}}(\cdot) + \nabla_{\bs{X}}(\cdot)^T \right )$. 
Mirroring the discretization strategy of section~\ref{meth:solid_mech_num}, we can derive the discrete weakform equation for the problem in Eq.~\eqref{eq_conservation_momentum_lin}-\eqref{eq_bc_lin} as: For every $ n = 1, \ldots N_T $, find $ (\bs{u}_h^n, \bs{v}_h^n, p_h^n) \in \bs{\mathcal{U}}^h \times \bs{\mathcal{V}}^h_0 \times \mathcal{P}^h $ such that,
\begin{equation} \label{eq_disc_wf_lin}
    \left ( \varrho \delta_t \bs{v}_h^n , \bs{y}_h \right )
    +
    ( \bs{\Sigma}^n, \nabla_{\bs{X}} \bs{y}_h )
     -
    (\bs{b}^n, \bs{y}_h)
    +
    (q_h,\nabla_{\bs{X}} \cdot \bs{v}_h^n )
    +
    (\bs{v}_h^n - \delta_t \bs{u}_h^n, \bs{w}_h)
    = 
    0
\end{equation}
for all $ (\bs{w}_h, \bs{y}_h, q_h) \in \bs{\mathcal{U}}^h \times \bs{\mathcal{V}}^h_0 \times \mathcal{P}^h $.
Here we assume $ \bs{u}_h^0 = \bs{u}_0$, $ \bs{v}_h^0 = \bs{v}_0 $, and $ \bs{Q}_k^0 = \bs{0}$ for each $ k = 1,\ldots N $ (\emph{e.g.} no history of deformation is assumed).
$ \bs{\Sigma}^n $ is given by Eq.~\eqref{eq_stress_lin}, which utilizes the approximate fractional derivative denoted by $ \hat{\text{D}}_n^\alpha $.  
\begin{equation}\label{eq_stress_lin}
    \bs{\Sigma}^n = E \bs{D} \bs{u}_h^n + \eta \hat{\text{D}}_n^\alpha (\bs{D} \bs{u}_h) + p_h^n \bs{I}
\end{equation}
Examining the stability of Eq.~\eqref{eq_disc_wf_lin}, it can be shown that the state variables satisfy the stability estimates given in Lemma~\ref{lem_lin_stability} and~\ref{lem_lin_stability_p} in~\ref{app:stability}.

\section{Approximation of Polynomial Functions}\label{sect:poly_example}

    In this section, we aim to examine the accuracy and convergence of the Prony-based method (Eq.~\eqref{eq_prony_approx_final}) in comparison to a few other methods in literature.
    For this purpose, we here consider a basic example of approximating the fractional derivative of a polynomial series, i.e.
\begin{equation}\label{eqn:polynomialfunction}
    p(t) = \sum_{k=1} b_k t^{k-1},
    \quad
    D_t^\alpha p(t) = \sum_k b_k \Gamma(1-k)\frac{k-1}{k-\alpha} t^{k-1-\alpha},
\end{equation}
    where $b_k$ are the given parameters (given under Fig.~\ref{fig:polyfit}). 
    Different methods are used to approximate the fraction derivative of $ p $, where the $L^2-\mathrm{norm}$ of the error with the analytic solution are then computed.
    All methods considered in this section were implemented in MATLAB  \cite{matlab2018version}. 
    
    Table~\ref{tab:polynomialtime} compares the computation time of the Prony-based method in Eq.~\eqref{eq_prony_approx_final} with a representative cumulative approach, in this case, the midpoint (MP) rule  (Eq.~\eqref{eqn:basiccaputo}).
    From this, we observe that the computation time for the midpoint rule grows as $ \mathcal{O}(N_T^2) $, with $ N_T $ denoting the number of time steps.  
    As a consequence, the midpoint rule Eq.~\eqref{eqn:basiccaputo} takes $ \approx 13$ min to complete the problem with $\Delta_t= $1e-5 in comparison to about $3$s for the Prony-based method. 
    This disparity stems from the use recursive updates of intermediate variables, $ q_k $, in the Prony-based method relies to make the temporal scaling $ \mathcal{O}(N_T) $.  
    Minimal increase in compute time is observed with increasing Prony terms despite the $N$-dependent workload required for the computation of the fractional derivative operation (computation time grows $15 \% $ from $ N = 3$ to $ N = 12 $).
    This likely stems from the efficiency of vectorized operations in MATLAB.

    For the Prony series parameters, we selected the optimized values from section \ref{meth:frac_opt} on a time scale of 9 seconds, 10 times the time domain. This can significantly improve the error estimates in Theorem \ref{thm_prony_error} by reducing $\epsilon(z)$ and $C(\bs{\beta},\bs{\tau})$, and being more inclusive of the frequencies of the Prony series that are larger than the time domain but still relevant to the problem (See Appendix \ref{app:timescaling} for more details).

    Fig.~\ref{fig:polyconverge} illustrates the accuracy of the proposed, midpoint and Gr{\"{u}}nwald-Letnikov methods. 
    Here the accuracy of the Prony approximation is examined for different $ \alpha $, $ N $ and $ \Delta_t $ ($\Delta_t = $1e-4 case shown in Fig.~\ref{fig:polyfit}). 
    Visibly, six or more Prony terms appears to perform well for the polynomial function, but three Prony terms generally yields noticeable errors. 
    As predicted by the truncation error in Fig.~\ref{fig:birkprony_err_comp}, the accuracy of the Prony approximation deteriorates with increasing $\alpha$, especially for $\alpha>0.5$, limiting the $ \Delta_t $ convergence as shown in Theorem~\ref{thm_prony_error}. 
    However, this deterioration is also observed in the midpoint method, where limited convergence with $\Delta_t$ is observed at $\alpha = 0.8$. The Gr{\"{u}}nwald-Letnikov method is significantly better in this regard, where convergence of the Gr{\"{u}}nwald-Letnikov method is on par or better than the Prony-based methods for $N=12$ for sufficiently large $N_T$. However, as observed in Fig.~\ref{fig:polyconverge}, the Gr{\"{u}}nwald-Letnikov method seems to diverge as $\Delta_t$ approach 0. This is due to a singular behavior in the Gr{\"{u}}nwald-Letnikov method (Eq.~\eqref{eq_gl}) when the number of time points $N_T$ is small. Precisely, when $N_T=0$, i.e. when taking the fractional derivative at time zero, the resulting value behaves as $f(0)/\Delta_t^\alpha$ (Fig. \ref{fig:polyfit}). This singularity is always present when $N_T$ is small and decreases quickly with increases $N_T$, but this also indicates the reliance of this method on knowing the value of the function being differentiated, making it unsuitable for solving differential equations without further corrections. 
    The Prony-based method for $N>3$ generally shows convergence for all $ \alpha $ values, with errors plateauing based on the truncation error derived in Lemma~\ref{lem_Prony_error} (see table~\ref{tab:polynomialconverge}).
    This is most visibly obvious for $\alpha=0.1$ (Fig. \ref{fig:polyconverge}A). 
    For the $N=3$ case, truncation errors vastly exceed time discretization errors, resulting in no appreciable convergence with time step refinement.
    It is noticeable that the Prony-based method performs better than the midpoint rule before reaching the truncation error in all cases tested. 
    Larger number of Prony terms can significantly reduce this truncation error, illustrating the importance of choosing a sufficient number of terms for your problem depending on the required degree of accuracy. 
    
    In addition to the midpoint rule, we also examined the method of Birk and Song~\cite{birk2010improved}, which also uses a Prony series approximation.
    The method presented in this work performs similarly or better for all $\alpha$ and $N$ (see Table~\ref{tab:polynomialconverge}). 
    This is consistent with our theoretical analysis in  Fig.~\ref{fig:birkprony_err_comp}, which shows that the method presented here has better performance with $ N > 3 $ in nearly all scenarios. 
    This improved performance stems from the increased error in the Prony approximation, causing errors to plateau with time step refinement earlier as predicted in  Lemma~\ref{lem_Prony_error} (see table~\ref{tab:polynomialconverge}).
    The performance gain is more significant for larger $\alpha$ and larger $N$ as predicted. 
    In the Prony-based method presented in this work, the integral is done with Prony terms with frequencies that are comparable with domain of the simulation, whereas in the method of Birk and Song approximates the complete range of frequencies by Gauss-Jacobi weights (which need not be optimal). 
    As a result, the method of Birk and Song does not see as much improvement in accuracy with an increases in $N$.

\section{FE approximation of fractional differential equations}
\label{sect:gao_example}

    While section~\ref{sect:poly_example} demonstrates the efficacy of the Prony approximation for given known differentiable functions, it does not give an indication of how well the method performs for solving a fractional differential equation.
    To examine this, we replicated Example 1 of Gao \textit{et al.} \cite{gao2012finite}, i.e.
\begin{equation} \label{eq:diffeqexample}
\begin{aligned}
    &D_t^\alpha u - \dpd[2]{u}{x} = f(x,t), \quad \mathrm{for} \, x\in[0,1], \quad t\in[0,1], \\
    &f(x,t) = (x-1)^2 e^{-x} t^3 \left[\Gamma(4+\alpha)\frac{(x-1)^2}{6} - (21-10x+x^2) t^\alpha\right] - 4x^2(x-1)^2(14x^2 - 14x+3),\\
    &u(x,0) = x^4 (x-1)^4, \\
    &u(0,t) = t^{3+\alpha}, \quad u(1,t) = 0.
\end{aligned}
\end{equation}
    The problem above has the unique analytic solution given in Eq.~\eqref{eq_fde_solution} and shown in Fig.~\ref{fig:FEfit}. 
    The major feature is a cubically growing source term at $x=0$, which decays to 0 at $x=1$ with a smaller concave quartic polynomial spanning $x\in[0,1]$. 
\begin{equation} \label{eq_fde_solution}
    u(x,t) = (x-1)^4 \left(e^{-x}t^{3+\alpha} + x^4\right) 
\end{equation}
    Following Gao \textit{et al.} \cite{gao2012finite}, Eq.~\eqref{eq:diffeqexample} is solved using linear finite elements in space.  
    Further, to make results comparable, we use the same error norm, i.e. 
\begin{equation}
    E(h, \Delta_t) = \max_{x \in [0,1]} \max_{t \in [0,1]} \left|u(x,t) - u_\mathrm{analytic}(x,t)\right| 
\end{equation}
    and replicate the results from Table 1 (temporal refinement) and Table 2 (spatial refinement) of Gao \textit{et al.} \cite{gao2012finite} (see Tables~\ref{tab:gaotable1} and~\ref{tab:gaotable2},  respectively). We also obtained Prony series parameters for the time scale of $[0, 100]$, by a factor of 100 times the time domain of the FE problem, to improve the error estimates (See Appendix \ref{app:timescaling}).
    Spatiotemporal errors are also shown in Fig.~\ref{fig:FEres} for $\Delta_t = 1/2500$, $h = 1/2500$ and a varying number of Prony terms.
    As expected, the error improves with an increasing number of Prony terms.

    Examining temporal refinement (with $h = $ 5e-5), the Prony-based method presented here has better convergence at smaller $\Delta_t$; however, the results from Gao \textit{et al.} show smaller initial errors (see Table~\ref{tab:gaotable1}). 
    This improved error for larger $ \Delta_t $ comes with additional computational expense requiring storage of the solution for all time points. 
    Additionally, we observe that the truncation error due to the Prony approximation results small errors, meaning that temporal convergence is observed even for $ N = 3$ up to $ \Delta_t > 1/160 $.
    Further, at coarse time discretizations, the Prony-based method shows increasing errors for increasing $ N $.
    While counter-intuitive, this is explained by Theorem~\ref{thm_prony_error}.
    With increasing $ N $, the discretization error is scaled by $C(\bs{\beta}, \bs{\tau})$ which generally grows with $ N $. This results in slow convergence until the approximation errors asymptotically converges to the error bounds. 
    However, optimal convergence is still achieved with $ \Delta_t $ refinement.

    In the case with refinement in space (with $\Delta_t =$ 5e-5), both methods have similar convergence rates but the Prony-based approach has smaller initial errors due to better convergence with $\Delta_t$ (Table~\ref{tab:gaotable2}). 
    Interestingly, but not unexpectedly, the number of Prony terms does not affect the convergence rate with refinement in space due to the fractional derivative being taken with respect to time only. 
    However, additional Prony terms are still important, as an insufficient number of terms limits the lower bound of the errors, as seen for $N=3$. 

    
    Fig.~\ref{fig:gao4}A, B and C shows spatiotemporal convergence for a range of time and discretizations.
    For fine spatial grids, we observe that errors approach the limit predicted by Lemma~\ref{lem_disc_Prony_error} and adhere to the results of Theorem~\ref{thm_prony_error}.
    Coarsening spatial refinement, we can observe spatial errors shift the lower bound for higher values of $ N $.
    Fig.~\ref{fig:gao4}D, E, and F presents a similar analysis utilizing the weights from Birk and Song~\cite{birk2010improved}.
    Once again, the Prony-based method presented in this work performs better, where the convergence rates are improved and lower bounds for the Prony-based approximation one ($N=9$) to two ($N=3$) orders of magnitude lower.

\section{Fractional Approximation applied to  Biomechanics in the Liver}
\label{sect:liver_example}
\subsection{Finite element model for the mechanical testing of livers}
In order to investigate the efficacy of the proposed method for simulating fractional viscoelasticity in soft tissues, here we model an experiment executed by Tan \emph{et al.}~\cite{Tan2013} on \textit{ex vivo} bovine liver tissue. 
Briefly, cylindrical liver samples of height $H = 2.7 $mm and radius $ R = 10 $mm were tested in a rheometer at a range of compressions, shear strains and frequencies. 
The behavior of the \emph{ex vivo} liver sample can be effectively modeled using a viscoelastic exponential model~\cite{Capilnasiu2019} of the form
\begin{equation}\label{eq_liver_pk2}
    \bs {S} = \delta \text{Dev} \left [D_t^\alpha \bs{S}_{v} \right] + p J \bs{C}^{-1}, \quad \bs{S}_{v} = \exp[b(II_{\bs{C}}-3)]\bs{C}.
\end{equation}
The first term denotes the deviatoric component of the stress depending on the fractional derivative of $ \bs{S}_v$ which follows a similar form to a typical hyperelastic material model of the right Cauchy-Green strain (see section~\ref{sect:kinematics}).
The second component denotes the hydrostatic stress, depending on the pressure variable, $ p $.
Here $\text{Dev}[\bs{A}] = \bs{A} - (\bs{A}:\bs{C} / 3) \bs{C}^{-1} $ denotes the deviatoric operator.
The model can be fit to match the behavior of the ex vivo liver from~\cite{Tan2013} by three parameters: the linear scaling parameter $ \delta $, the nonlinear response parameter $ b $, and the fractional derivative $ \alpha $. 
Following the parameterization process in~\cite{Capilnasiu2019}, parameters $\alpha=0.2$, $b=1.5$ and $\delta=126.4$Pa were observed to provide the best fit to the data (see Fig.~\ref{fig:liver_behaviour}).

Here we focus on a single test case, where a sample was compressed by $ 10 \% $ and subsequently sheared at a frequency of $ 1 $Hz around its axis with a shear strain amplitude of $ 25 \% $.
In this case, the reference domain is the cylindrical region $ \Omega_0 = \{ \bs{x} \in \mathbb{R}^3 | \;  \sqrt{x_1^2 + x_2^2} \le R, 0 \le x_3 \le H \} $. 
The boundary $ \Gamma_0 $ is split into a bottom surface, top surface, and side surface $ \Gamma_0 = \Gamma_0^b \cup \Gamma_0^t \cup \Gamma_0^w $, respectively.
Assuming pure compression and shear, the deformation can be described in analytic form as,
\begin{equation}\label{eq_liver_analytic_sol}
	\bs u (t)=
	\left[
	\begin{array}{c}
	 \dfrac{X_1 }{\sqrt {\lambda(t)}} \cos (\psi(t) X_3)   -\dfrac{ X_2 }{\sqrt {\lambda(t)}} \sin (\psi(t) X_3) - X_1 \vspace{2mm} \\
	  \dfrac{X_1}{\sqrt {\lambda(t)}} \sin (\psi(t) X_3)  + \dfrac{ X_2 }{\sqrt {\lambda(t)}} \cos (\psi(t) X_3)  - X_2   \vspace{2mm}\\
	X_3 \left ( \lambda(t)-1 \right )\\
	\end{array}
	\right],
	\quad
	\psi(t)= \frac{\gamma \lambda^{3/2}(t) \sin(2\pi f \hat{t})   }{\sqrt{X_1^2 + X_2^2}},
\end{equation}
where the compression, $ \lambda $, and the shear angle, $ \psi X_3 $, were chosen to mimic the experimental conditions.
Here compression took place over the first $1$s, \emph{e.g.} $\lambda (t) = 1-0.1\min(t,1) $. 
Additionally, the angular displacement was applied after compression, $ \hat{t} = \max(0,t-1) $, with an oscillatory frequency $f = 1$Hz and a maximal shear strain of $\lambda \gamma$.
Torque readings were acquired during the shearing process. 

This problem provides the ability to verify the finite element implementation outlined in section~\ref{meth:solid_mech_num}.
The analytical stresses could be computed using Eq.~\eqref{eq_liver_analytic_sol}) and Eq.~\eqref{eq_liver_pk2}.
Applying Eq.~\eqref{eq_cauchystress}, the model predicted torque, $ \tau $, and normal traction, $ t_N $, can be calculated by
\begin{equation}
    \tau (t) =  \int_{\Gamma^t} \bs{e}_3 \cdot \left [ \bs{x} \times (\bs{\sigma} \cdot \bs{n}) \right ] \; d \Gamma,
    \quad
    t_N (t) = \int_{\Gamma^t} \bs{e}_3 \cdot \left [ \bs{\sigma} \cdot \bs{n} \right ] \; d \Gamma ,
\end{equation}
where $ \Gamma^t $ is upper boundary of the physical domain $ \Omega $, $ \bs{n} = \bs{e}_3 $ the outward normal, and $ \bs{e}_k $ is the $k^{th}$-base vector for rectangular-Cartesian coordinates.
These quantities were computed in MATLAB R2018b and in-house software $\mathcal{C}$\textbf{heart}. 
In MATLAB, a reference solution for the torque and normal traction were computed using the Gr{\"{u}}nwald-Letnikov method~\cite{Scherer2011} (Eq.~\eqref{eq_gl}) with a small time-step $\Delta_t=$1e-5.
A sweep over time-step size and the number of Prony terms was also applied to examine computing times and computational accuracy.

For the finite element implementation in $\mathcal{C}$\textbf{heart}, a curvilinear hexahedral mesh (with $9,072$ elements) was used to approximate $ \Omega_0 $.
The inf-sup stable $ \mathbb{Q}^2 - \mathbb{Q}^2 - \mathbb{Q}^1$ scheme was used to solve for $ \bs{u} $, $ \bs{v} $ and $ p $ (with $77,425$ quadratic nodes and $ 10,309 $ linear nodes).
Eq.~\eqref{eq_liver_analytic_sol} was applied as boundaries conditions on $\Gamma_0^b$ and $\Gamma_0^t$ with no boundary conditions applied on $ \Gamma_0^w$.  
Model torque and normal traction were exported along with the solution.

\subsection{Accuracy and convergence}
The accuracy and convergence of the Gr{\"{u}}nwald-Letnikov and Prony approximations are shown in Table~\ref{tab:newaccuracy}, with the latter examined for $N \in \{ 3, 6, 9, 12 \} $.
When the number of Prony terms $ N = 3 $, there is little improvement in the stress accuracy with temporal refinement reflecting inherent limitations in the approximation quality with so few terms.
However, increasing to $ N = 6$, we see the stagnation in convergence shift to $ \Delta_t =$1e-3.
Further increases to $ N = 9 $ and $ N = 12 $ exhibit continued convergence even for the finest time step size, s predicted from Theorem \ref{thm_prony_error}

The deformation and pressure from the FEM simulation can be observed in Fig.~\ref{fig:liver} (with $N = 9$ and $ \Delta_t =$1e-3).
The simulated displacements were within 1e-6 of their analytic counterparts. The FEM simulation was further verified against the Gr{\"{u}}nwald-Letnikov implementation in MATLAB.
The torque and normal force from both implementations are shown in Fig.~\ref{fig:liver}, showing both forward model and the analytically computed behavior are consistent.  Here, an error of 1e-3 was observed for the torque and 2e-4 was observed in the normal traction force.

\subsection{Computational cost}

Table~\ref{tab:oldnewtimecomparison} presents the cost performance of the Prony approximation (Eq.~\eqref{eq_prony_approx_final}) compared to the Gr{\"{u}}nwald-Letnikov method (Eq.~\eqref{eq_gl}) over the time domain $t \in [0,2]$. 
Time step $\Delta_t$ was chosen to vary between 1e-2 and 5e-5. 
For time steps larger than 1e-4, the the Gr{\"{u}}nwald-Letnikov and Prony approach (for $n=3, 6, 9, 12$ terms) are comparably fast. 
However, as $ \Delta_t $ shrinks further, the Gr{\"{u}}nwald-Letnikov  approach becomes increasingly expensive to compute ($ \mathcal{O} (N_T^2)$ versus $ \mathcal{O}(N_T) $), with storage costs of $ \mathcal{O}(N_T) $ vs $ \mathcal{O}(N) $ (N is the number of Prony terms).

Finally, we sought to examine the added computational cost of simulating the viscoelastic model compared to a classical hyperelastic model (see Table~\ref{tab:speedelasticviscoelastic}).
Here the compute time of the viscoelastic model in Eq.~\eqref{eq_liver_pk2} was contrasted against the hyperelastic model resulting if $ \alpha = 0 $ in Eq.~\eqref{eq_liver_pk2}.
It is notable that the computational time for both hyperelastic and viscoelastic material models does not increase with direct proportionality with the number of time steps, likely due to a reduction in the required number of Newton-Raphson iterations per time step.
Further, while the computational cost of the viscoelastic model is larger, this does not lead to a significant increase in computational time, with the differences observed being within $\approx$16\%.
Additionally, varying the number of Prony terms does not significantly impact the computational time, as seen in Table~\ref{tab:speedelasticviscoelastic}.
While increasing $ N $ increases the computational work, the added cost for residual and matrix calculations is minor relative to the basic computational expense for simulating large deformation mechanical systems.


\section{Discussion}
\subsection{Computational cost}
    The biggest advantage of this Prony based method for approximating Caputo fractional derivatives Eq. \eqref{eq_prony_approx_final} is computational efficiency. The historical integral required typically leads to an accumulating computational cost and demands smaller time steps by introduces nonlinearity into the system in the form of the strain history, becoming a heavy burden in viscoelastic simulations. The method developed herein allows a bounded storage cost for the memory cost of the strain history by an approximation using Prony series. The parameters can be precomputed from the $\alpha$ values and are independent of the strain history itself, which only leaves the cost of summing over the Prony series at each time step, which is small. The number of Prony terms necessary is dependent on the required accuracy of the approximation. However, for the example we have tested, six or fewer terms are sufficient for a wide range of cases. We find this cost to be only a fraction of the cost of the overall simulation, but can reduce the cost of simulating later time steps by several magnitudes. A well-bounded cost for the historical integral also removes a limitation on the full duration of the simulation, making it applicable to a wider range of applications.

\subsection{Accuracy and convergence}
    The number of terms to use in the Prony-based method (Eq. \eqref{eq_prony_approx_final}) is the decision between the truncation error, lower-bound for the approximation, and the loss in accuracy at larger $\Delta_t$ due to the $C(\bs{\beta}, \bs{\tau})$ term. Increase in the number of Prony terms with lower the resulting truncation error at the cost of larger error terms. The Prony-based method has superior computational costs and better convergence rates when above the truncation error (Lemma~\ref{lem_Prony_error}) in comparison to classical non-recursive methods (Fig.~\ref{fig:polyconverge}). It is noted that in some cases, the discretization error introduced by the coefficient $C(\bs{\beta},\bs{\tau})$ (Lemma~\ref{lem_disc_Prony_error}) can cause the solution to be initially be less accurate at large $\Delta_t$ (Table~\ref{tab:gaotable1} and \ref{tab:gaotable2}). However, this factor remains to be small and the savings in the computational cost is the much more significant factor. Against competing methods such as the Birk and Song method \cite{birk2010improved}, the Prony-based method has better convergence 
    for six or more Prony terms or quadrature points and at higher $\alpha$, 
    and is competitive otherwise (Fig.~\ref{fig:birkprony_err_comp}) in examples tested.

\subsection{Other methods}
    
    One common approach to improve the computational efficiency of traditional methods is by truncating the historical data. As noted by Diethelm et al. \cite{diethelm2005algorithms}, it is typical that the weights in fractional derivative approximations exhibit some form of fading memory. For this reason, terms in the history are neglected to reduce computational costs. While it is true that terms that are further in the past will decay more, and in a way becomes small enough to be negligible, the length of the history also increases. As such, truncation cannot be done negligently. For example, for $\alpha=0.3$ and if the method presented by Diethelm et al. \cite{diethelm2005algorithms} in Eq. \eqref{eqn:diethelgrunwald} is truncated to 12 terms, the sum of the absolute value of the first 12 weights on the domain $t\in[0,1]$ is always 1.66353, independent of $\Delta_t$ (assuming $\Delta_t < 1/11$). In comparison, the sum of the remaining weights are 0.16037, 0.248329, 0.292299, and 0.314331 for $\Delta_t = $ 1e-2, 1e-3, 1e-4 and 1e-5 respectively. The weight of the truncated history actually increases faster than the rate of decay due to fading memories. This effect is more significant for lower values of $\alpha$, and thus care must be taken for how the truncation is done. The number of terms to keep need to increase over time, depending on the problem.
    
    Another similar method has been presented by Jiang et al. \cite{jiang2017fast}. The major differences between the two methods are: i) the approximation of the integral in the discretized case and ii) the method of computing the parameters of the Prony series. For the approximation of the integral, this is largely a choice. However, we do note that in the method of Jiang \textit{et al.} \cite{jiang2017fast}, by holding the input function constant will still lead to a small but non-zero increase in the history terms when the solution should be zero on a stationary pure viscoelastic body. Jiang \textit{et al.}'s method of computing the Prony series parameters \cite{jiang2017fast} is designed with the goal of targeting a desired error of the approximation. However, the bounds on the number of terms are typically in the hundred for even modest desired error (e.g., 1e-4). 
    The algorithm provided by Jiang \textit{et al.} \cite{xu2012bootstrap} can reduce this number, but it is still larger than the Prony-based method. In addition, the algorithm presented can select from a list of frequencies based on Eigenvalues of the least squares system, which do not necessarily include the optimal frequencies.
    The method presented in this work (Eq. \eqref{eq_prony_approx_final}) outputs a desired number of Prony terms, which is more desirable for ease of implementation and application.

\subsection{Other applications}

    Although this method is described for the constitutive modeling of viscoelastic materials, fractional derivatives have a wide range of other applications. The most common application of fractional calculus in the literature is in regard to viscous fluid flows. The parallels to viscoelasticity is obvious and its use along with viscoelasticity to solve complex biomechanical systems is also very common. Also common is the application of acoustic waves traveling through viscous media. Another interesting and useful application is in edge detection in image analysis. The use of fractional derivatives can improve thin edge detection and tolorence to noises \cite{mathieu2003fractional}. One other important application is in PID controllers, where the feedback loops are on fractional order. An interesting example of this is in path tracking of self driving vehicles \cite{suarez2003using}. Some other common and interesting applications are summarized by Dalir and Bashour \cite{dalir2010applications}.

\subsection{Practical guidance}
    
    The cost of computing the parameters of the Prony series can be expensive and are best precomputed then imported at the time of use. However, since each of the parameters as a function of $\alpha$ is smooth, one method is to generate a set of Prony series parameters for well-spaced values of $\alpha$, which can then be interpolated for more specific values of $\alpha$. Indeed, we used this method for the MATLAB and Cheart \cite{lee2016multiphysics} implementations. The number of Prony terms to use is quite flexible. $N<4$ is not recommended, but six is typically sufficient for most problems. Increasing $N$ further will also not significantly negatively impact the solution and the computational cost of the problem. It should be noted that due to the necessity of approximating the derivative of the function to be uniform for each time interval, smaller time steps are necessary for the convergence of the solution in finite element simulations. One should choose step sizes such that the derivative is expected to be constant. Also, we have found that the velocity is unstable under second-order schemes in transient viscoelasticity simulations, but is stable under first-order schemes, as shown in \ref{app:stability}.
\section{Conclusion}

In this work we present a novel method for the finite element approximation of fractional nonlinear viscoelastic materials.
Application of fractional viscoelastic models can present computational challenges, with many methods requiring prohibitive storage and computational expense.
To overcome these issues, a new approximation to the Caputo fractional derivative is introduced, resulting in a fixed storage cost and making the computational times scale linearly with the number of time steps.
Error estimates are derived for this approximation, demonstrating the efficacy of the approximation.
Estimates are demonstrated to be accurate through numerical examples, where comparable or better convergence behaviors with refinement in time are observed compared to current state-of-the-art methods. 
The fractional approximation was then introduced into a nonlinear transient finite element mechanics implementation.
The numerical problem was shown to be unconditionally stable for simple linear fractional viscoelasticity, and observed stable for nonlinear viscoleasticity through verified numerical examples.
The computational cost of this fractional derivative implementation is comparable to the standard hyperelastic model, despite the historical integral ($\sim$16\% more), and it was able to accurately approximate the solution in all examples considered.
This method can significantly simplify the implementation of fractional viscoelastic models in general, and provides a pathway for applications of viscoelastic soft tissue models for the efficient simulation of complex biomechanical systems.
While developed in the context of viscoelastic modelling, the Caputo approximation introduced has potential applications in solving the many other problems involving fractional derivatives. 




\section{Acknowledgements}
\begin{sloppypar}
\noindent Authors would like to acknowledge funding from Engineering and Physical Sciences Research Council (EP/N011554/1 and EP/R003866/1), and support from the Wellcome Trust EPSRC Centre of Excellence in Medical Engineering (WT 088641/Z/09/Z) and the NIHR Biomedical Research Centre at Guy's and St.Thomas' NHS Foundation Trust and KCL. 
The views expressed are those of the authors and not necessarily those of the NHS, the NIHR, or the DoH.  
\end{sloppypar} 

\section*{References}
\bibliographystyle{elsart-num-sort}
\bibliography{ref}

\newpage
\appendix 
\section{Stability Analysis for Linear Fractional Viscoelastic Materials} \label{app:stability}

Lemma~\ref{lem_lin_stability} and~\ref{lem_lin_stability_p} provide stability estimates for the linear viscoelastic mechanics problem introduced in Eq.~\eqref{eq_disc_wf_lin}.
We note that the fractional derivative approximation follows the form discussed in section~\ref{meth:frac_approx_disc}.
Here $ N $ denotes the number of Prony terms, and $ \beta_k, \tau_k \in \mathbb{R}^+ $ denote the scaling and time scale of the $k^{th} $ Maxwell element in the approximation (with $ \beta_0 $ being the scaling of a pure dashpot element).
This approximation relies on intermediate variables $ \bs{Q}_k^n $ which obey the update formula in Eq.~\eqref{eq_stress_frac_lin}.
\begin{equation}\label{eq_stress_frac_lin}
    \hat{\text{D}}_n^\alpha (\bs{D} \bs{u}_h) 
    = \beta_0 \bs{D} \bs{v}_h^n + \sum_{k=1}^N \bs{Q}_k^n, 
    \quad
    \bs{Q}_k^n = e_k^2 \bs{Q}_k^{n-1} + \beta_k e_k \Delta_t \bs{D} \bs{v}_h^n 
\end{equation}
With this approximation in mind, we can derive the following stability estimate for the discrete solution showing that the discrete solution is unconditionally stable and bounded by given data.

\begin{lem}\label{lem_lin_stability}
Consider the linear incompressible fractional viscoelasticity problem shown in Eq.~\eqref{eq_disc_wf_lin} and the update formula shown in Eq.~\eqref{eq_stress_frac_lin}.  
Assuming that $ \bs{v}_h^0 , \bs{u}_h^0 \in \bs{\mathcal{V}}^h_0  $,
that $ \bs{Q}_k^0 = \bs{0}$,
and 
$$ \bs{b} \in \bs{L}^{\infty} ([0,T]; \bs{L}^2(\Omega_0)), \quad \sup_{t \in [0,T]} \| \bs{b}(t) \|_0 \le K_{\bs{b}} $$
%
\begin{equation}
    \varrho \| \bs{v}_h^n \|_0^2      
    +
    E \| \bs{D} \bs{u}_h^n \|_0^2
    +
    \eta \beta_0 \sum_{m=1}^n \Delta_t \| \bs{D} \bs{v}_h^m \|_0^2
    +
    \sum_{k=1}^N  \frac{\eta}{\beta_k e_k}  \| \bs{Q}_k^n \|_0^2  
    \le 
    \varrho \| \bs{v}_h^{0} \|_0^2
    +
    E \| \bs{D} \bs{u}_h^{0} \|_0^2
    +
    \frac{C_\Omega t_n K_{\bs{b}} }{\eta \beta_0} 
    \nonumber
\end{equation}
for each $ n = 1, \ldots N_T$.
\end{lem}

\begin{pol}{\ref{lem_lin_stability}}
To prove Lemma~\ref{lem_lin_stability}, we choose $ \bs{w}_h = \bs{0} $, $ \bs{y}_h = \bs{v}_h^n $ and $ q_h = 0 $ in Eq.~\eqref{eq_disc_wf_lin} resulting in the equation
\begin{equation} \label{eq_lem_stab_1}
    \left (\varrho \delta_t \bs{v}_h^n , \bs{v}_h^n \right )
    +
    ( \bs{\Sigma}^n, \nabla_{\bs{X}} \bs{v}_h^n )
    -
    (\bs{b}^n, \bs{v}_h^n)
    = 
    0
\end{equation}
Rearranging the update formula for $ \bs{Q}_k^n $ in Eq.~\eqref{eq_stress_frac_lin}, we note that
\begin{equation} \label{eq_lem_stab_2}
    \bs{D} \bs{v}_h^n = \frac{1}{\beta_k e_k \Delta_t} \left ( \bs{Q}_k^n - e_k^2 \bs{Q}_k^{n-1} \right).
\end{equation}
Focusing on the stress term in Eq.~\eqref{eq_lem_stab_1}, noting the symmetry of $ \bs{\Sigma}^n $,that $ \bs{v}_h^n $ is weakly divergence free, the modified update formula of Eq.~\eqref{eq_lem_stab_2}, and the identity
\begin{eqnarray} 
    (\bs{a} - \bs{b}) \cdot \bs{a} 
    = 
    \frac{1}{2} \left( | \bs{a} |^2 - | \bs{b} |^2 + | \bs{a}-\bs{b} |^2 \right),
    & \; &
    \text{ for any }
    \bs{a},\bs{b} \in \mathbb{R}^d
    \label{eq_lem_stab_3} \\
    (\bs{A} - \bs{B}) : \bs{A} 
    = 
    \frac{1}{2} \left( | \bs{A} |^2 - | \bs{B} |^2 + | \bs{A}-\bs{B} |^2 \right),
    & \; &
    \text{ for any }
    \bs{A},\bs{B} \in \mathbb{R}^{d \times d}
    \nonumber
\end{eqnarray}
we can observe that
\begin{eqnarray}
    ( \bs{\Sigma}^n, \nabla_{\bs{X}} \bs{v}_h^n )
    & = &
    ( \bs{\Sigma}^n, \bs{D} \bs{v}_h^n )
    \nonumber \\
    & = &
    \left ( E \bs{D} \bs{u}_h^n + \eta \hat{\text{D}}_n^\alpha (\bs{D} \bs{u}_h) + p_h^n \bs{I} , \bs{D} \bs{v}_h^n \right )
    \nonumber \\
    & = &
    \frac{E}{\Delta_t} \left( \bs{D} \bs{u}_h^n , \bs{D} [\bs{u}_h^n - \bs{u}_h^{n-1}] \right)
       +
       \eta \left( \beta_0 \bs{D} \bs{v}_h^n + \sum_{k=1}^N \bs{Q}_k^n , \bs{D} \bs{v}_h^n \right)
    \nonumber \\
    & = &
    \frac{E}{\Delta_t} \left( \bs{D} \bs{u}_h^n , \bs{D} [\bs{u}_h^n - \bs{u}_h^{n-1}] \right)
       +
       \eta \beta_0 \| \bs{D} \bs{v}_h^n \|_0^2
       +
       \sum_{k=1}^N  \frac{\eta}{\beta_k e_k \Delta_t} \left( \bs{Q}_k^n , \bs{Q}_k^n - e_k^2 \bs{Q}_k^{n-1} \right)
    \nonumber \\
    & = &
    \frac{E}{2 \Delta_t} \left( 
          \| \bs{D} \bs{u}_h^n \|_0^2  
          - 
          \| \bs{D} \bs{u}_h^{n-1} \|_0^2 
          +
          \| \bs{D} (\bs{u}_h^n - \bs{u}_h^{n-1}) \|_0^2
       \right)
       +
       \eta \beta_0 \| \bs{D} \bs{v}_h^n \|_0^2
       \nonumber \\
    &  & \hspace{5mm}
       +
       \sum_{k=1}^N  \frac{\eta}{2 \beta_k e_k \Delta_t} 
       \left(  
          \| \bs{Q}_k^n \|_0^2  
          -  
          \|    e_k^2 \bs{Q}_k^{n-1} \|_0^2 
          +
          \|   \bs{Q}_k^n -  e_k^2 \bs{Q}_k^{n-1} \|_0^2 
       \right)
       \label{eq_lem_stab_4}
\end{eqnarray}
Applying Eq.~\eqref{eq_lem_stab_3} and the equality in Eq.~\eqref{eq_lem_stab_4} to Eq.~\eqref{eq_lem_stab_1}, 
\begin{eqnarray}
    & &
    \frac{\varrho}{2} \left( 
        \| \bs{v}_h^n \|_0^2
        -
        \| \bs{v}_h^{n-1} \|_0^2
        +
        \| \bs{v}_h^n - \bs{v}_h^{n-1} \|_0^2
    \right)
    +
    \frac{E}{2} \left(
        \| \bs{D} \bs{u}_h^n \|_0^2  
        - 
        \| \bs{D} \bs{u}_h^{n-1} \|_0^2 
        +
        \| \bs{D} (\bs{u}_h^n - \bs{u}_h^{n-1}) \|_0^2
    \right)
    \hspace{30mm}
    \nonumber \\
    & & \hspace{10mm}
    +
    \eta \beta_0 \Delta_t \| \bs{D} \bs{v}_h^n \|_0^2
    +
    \sum_{k=1}^N  \frac{\eta}{2 \beta_k e_k} \left(  
        \| \bs{Q}_k^n \|_0^2  
        -  
        \|    e_k^2 \bs{Q}_k^{n-1} \|_0^2 
        +
        \|   \bs{Q}_k^n -  e_k^2 \bs{Q}_k^{n-1} \|_0^2 
    \right)    
    -
    \Delta_t \left ( \bs{b}^n , \bs{v}_h^n \right ) 
    = 
    0. 
    \label{eq_lem_stab_5}
\end{eqnarray}
%
Examining the final term in Eq.~\eqref{eq_lem_stab_5}, noting through both Korn and Poincar\'{e} inequalities there exists a $ C_\Omega > 0 $ such that,
\begin{equation*} 
  \| \bs{v}_h^n \|_0 \le C_\Omega \| \bs{D} \bs{v}_h^n \|_0
\end{equation*}
and applying Young's inequality (with $ \epsilon = \eta \beta_0 / C_\Omega $), then
\begin{eqnarray}
    (\bs{b}^n, \bs{v}_h^n)
    & \le &
    \frac{1}{2 \epsilon} \| \bs{b}^n \|_0^2 + \frac{\epsilon}{2} \| \bs{v}_h^n \|_0^2 
    \nonumber \\
    & \le &
    \frac{C_\Omega}{2 \eta \beta_0} \| \bs{b}^n \|_0^2 
    + 
    \frac{\eta \beta_0 }{2} \| \bs{D} \bs{v}_h^n \|_0^2 
    \label{eq_lem_stab_6}
\end{eqnarray}
Combining Eq.~\eqref{eq_lem_stab_5} and~\eqref{eq_lem_stab_6} and re-arranging terms, we observe that
\begin{eqnarray}
    & &
    \varrho \| \bs{v}_h^n \|_0^2      
    +
    E \| \bs{D} \bs{u}_h^n \|_0^2
    +
    \eta \beta_0 \Delta_t \| \bs{D} \bs{v}_h^n \|_0^2
    +
    \sum_{k=1}^N  \frac{\eta}{\beta_k e_k}  \| \bs{Q}_k^n \|_0^2  
    \nonumber \\
    & & \hspace{20mm}
    \le 
    \varrho \| \bs{v}_h^{n-1} \|_0^2
    +
    E \| \bs{D} \bs{u}_h^{n-1} \|_0^2
    +
    \sum_{k=1}^N  \frac{\eta}{\beta_k e_k} \| e_k^2 \bs{Q}_k^{n-1} \|_0^2 
    +
    \Delta_t \frac{C_\Omega}{\eta \beta_0} \| \bs{b}^n \|_0^2 
    \nonumber
\end{eqnarray}
Noting that $ e_k \le 1 $, applying induction, noting $ \bs{Q}_k^0 = \bs{0}$ and utilizing the boundedness of $ \bs{b} $, we arrive at the stability estimate. 
\end{pol}

To show stability for the pressure, $ p_h^n $, we assume that the spaces $ \bs{\mathcal{V}}^h_0 $ and $ \mathcal{P}^h $ are inf-sup stable, satisfying the condition ($ \beta > 0$) ~\cite{quarteroni2008numerical},
\begin{equation} \label{eq_infsup}
 \forall q_h \in \mathcal{P}^h \; \; \exists \; \bs{w}_h \in \bs{\mathcal{V}}_0^h, \; \bs{w}_h \neq \bs{0} : \quad ( q_h , \nabla_{\bs{X}} \cdot \bs{w}_h ) \ge \beta \| q_h \|_0 \| \bs{w}_h \|_1.
\end{equation}
In addition, we assume there is a $ \bs{v}_h^{-1}  \in \bs{\mathcal{V}}^h_0 $ satisfying
\begin{equation} \label{eq_vel_init}
 \| \Delta_t^{-1} (\bs{v}_h^{0} - \bs{v}_h^{-1}) \| \le K
\end{equation}
for some $ K \ge 0 $ (independent of $ \Delta_t $) and that the initial conditions satisfy the discrete weakform (with $ \bs{Q}_k^{-1} = \bs{0} $ for $ k = 1, \ldots N $), 
\begin{equation} \label{eq_disc_wf_lin_init}
    \left (\varrho \delta_t \bs{v}_h^0 , \bs{y}_h \right )
    +
    ( \bs{\Sigma}^0, \nabla_{\bs{X}} \bs{y}_h )
    -
    (\bs{b}^0, \bs{y}_h)
    = 
    0
\end{equation}
for all weakly divergence free $ \bs{w}_h \in \bs{\mathcal{V}}^h_{0,Div} $, with the space
$$ 
  \bs{\mathcal{V}}^h_{0,Div} = \left \{ 
     \bs{w}_h \in \bs{\mathcal{V}}^h_{0} 
     \; \big | \;
     (q_h,\nabla_{\bs{X}} \cdot \bs{v}_h^n ) = 0,
     \quad
     \forall 
     q_h \in \mathcal{P}^h
  \right \}.
$$
With these assumptions, we can prove unconditional stability for the discrete model pressure.

\begin{lem}\label{lem_lin_stability_p}
Suppose the assumptions of Lemma~\ref{lem_lin_stability} hold.
Assuming there is a $ \bs{v}_h^{-1}  \in \bs{\mathcal{V}}^h_0 $ satisfying Eq.~\eqref{eq_vel_init},
that $ \bs{Q}_k^{-1} = \bs{0} $, that 
$$ \bs{b} \in \bs{W}^{1,\infty} ([0,T]; \bs{L}^2(\Omega_0)), \quad  
   \sup_{t \in [0,T]} \| \bs{b} (t) \|_0 + \sup_{t \in [0,T]} \| \partial_t \bs{b} \|_0 \le K_{\bs{b}}^\prime, $$
and that the spaces $ \bs{\mathcal{V}}_0^h $ and $ \mathcal{P}^h $ are inf-sup stable,
then $ p_h^n $ satisfies the stability estimate,
\begin{eqnarray}
    \| p_h^n \|_0      
    \le
    \frac{C}{\beta} \left(  
        K
        +
        K_{\bs{b}}
        +
        t_n^{1/2} K_{\bs{b}}^\prime
        +
        \| \bs{v}_h^{0} \|_0
        +
        \| \bs{D} \bs{u}_h^{0} \|_0
        +
        \| \bs{D} \bs{v}_h^{0} \|_0
    \right)
    \nonumber
\end{eqnarray}
for a constant $ C > 0 $ (independent of $ h $ and $ \delta_t $) and for each $ n = 1, \ldots N_T$.
\end{lem}

\begin{pol}{\ref{lem_lin_stability_p}}

Looking at Eq.~\eqref{eq_disc_wf_lin}, we first separate components of $ \bs{\Sigma}^n $, choose $ \bs{w}_h = \bs{0} $ and $ q_h = 0 $ and rearrange the equation.  
Applying Cauchy-Schwartz inequality, we can arrive at the following inequality.
\begin{eqnarray} 
    ( p_h^n , \nabla_{\bs{X}} \cdot \bs{y}_h)
    & = &
    (\bs{b}^n, \bs{w}_h)
    -
    \varrho \left ( \delta_t \bs{v}_h^n , \bs{y}_h \right )
    -
    ( E \bs{D} \bs{u}_h^n + \eta \hat{\text{D}}_n^\alpha (\bs{D} \bs{u}_h) , \nabla_{\bs{X}} \bs{y}_h )
    \nonumber \\
    & \le &
    \| \bs{b}^n \|_0 \| \bs{y}_h \|_0
    +
    \varrho \| \delta_t \bs{v}_h^n \|_0 \| \bs{y}_h \|_0
    +
    \| E \bs{D} \bs{u}_h^n + \eta \hat{\text{D}}_n^\alpha (\bs{D} \bs{u}_h) \|_0 \| \nabla_{\bs{X}} \bs{y}_h \|_0
    \nonumber \\
    & \le & \left (
    \| \bs{b}^n \|_0 
    +
    \varrho \| \delta_t \bs{v}_h^n \|_0 
    +
    E \| \bs{D} \bs{u}_h^n \|_0 
    +
    \eta  \beta_0 \|  \bs{D} \bs{v}_h^n \|_0  
    + 
    \sum_{k=1}^N \eta \| \bs{Q}_k^n \|_0 
    \right ) \| \bs{y}_h \|_1
    \label{eq_lem_stab_p1} 
\end{eqnarray}
Noting that the spaces of discrete solutions satisfy the inf-sup condition, the inequality in Eq.~\eqref{eq_lem_stab_p1} can be simplified to provide an upper bound on the pressure, $ p_h^n $.
\begin{equation}
    \| p_h^n \|_0 
    \le
    \frac{1}{\beta}
    \left (
        \| \bs{b}^n \|_0 
        +
        \varrho \| \delta_t \bs{v}_h^n \|_0 
        +
        E \| \bs{D} \bs{u}_h^n \|_0 
        +
        \eta  \beta_0 \|  \bs{D} \bs{v}_h^n \|_0  
        + 
        \sum_{k=1}^N \eta \| \bs{Q}_k^n \|_0 
    \right)
    \label{eq_lem_stab_p2}
\end{equation}
On the RHS of Eq.~\eqref{eq_lem_stab_p2}, the terms involving $ \bs{b} $ remain bounded due to the implicit assumption that $ \bs{b} \in \bs{W}^{1,\infty} ([0,T]; \bs{L}^2(\Omega_0)) $.
Moreover, from Lemma~\ref{lem_lin_stability}, we know there exists a $ C_0 > 0$ (independent of $ h $ and $ \Delta_t $) such that,
\begin{equation}
    E \| \bs{D} \bs{u}_h^n \|_0  + \sum_{k=1}^N \eta \| \bs{Q}_k^n \|_0 
    \le
    C_0 \left(
        \| \bs{v}_h^{0} \|_0
        +
        \| \bs{D} \bs{u}_h^{0} \|_0
        +
        t_n^{1/2} K_{\bs{b}}
    \right).
    \label{eq_lem_stab_p3}
\end{equation}
The stability estimate in Eq.~\eqref{eq_lem_stab_p2} relies on the boundedness of the discrete time derivative of $ \bs{v}_h^n $ as well as the symmetric derivative of $ \bs{v}_h^n $.
To ensure these quantities remain bounded independent of $ \Delta_t $, we return to Eq.~\eqref{eq_disc_wf_lin}.
Subtracting the equation for $ t_{n-1} $ from that for $ t_n $, letting $ \delta (\cdot)^n =  (\cdot)^n  -  (\cdot)^{n-1} $, and selecting the test functions $ \bs{w}_h=\bs{0} $, $ \bs{y}_h = \delta \bs{v}_h^n $ and $ q_h = 0 $ yields
\begin{equation} 
    \left (\varrho \delta ( \delta_t \bs{v}_h^n ), \delta \bs{v}_h^n \right )
    +
    ( \delta \bs{\Sigma}^n , \nabla_{\bs{X}} \delta \bs{v}_h^n )
    -
    (\delta \bs{b}^n, \delta \bs{v}_h^n )
    = 
    0
    \label{eq_lem_stab_p4}
\end{equation}
Applying the same approach as in Lemma~\ref{lem_lin_stability} we can arrive at the following bound.
\begin{eqnarray}
    & &
    \varrho \| \delta \bs{v}_h^n \|_0^2      
    +
    E \| \bs{D} (\delta \bs{u}_h^n) \|_0^2
    +
    \eta \beta_0 \Delta_t \| \bs{D} (\delta \bs{v}_h^n) \|_0^2
    +
    \sum_{k=1}^N  \frac{\eta}{\beta_k e_k}  \| \delta \bs{Q}_k^n \|_0^2  
    \nonumber \\
    & & \hspace{20mm}
    \le 
    \varrho \| \delta \bs{v}_h^{n-1} \|_0^2
    +
    E \| \bs{D} (\delta \bs{u}_h^{n-1}) \|_0^2
    +
    \sum_{k=1}^N  \frac{\eta}{\beta_k e_k}  \| e_k^2 \delta \bs{Q}_k^{n-1} \|_0^2 
    +
    \Delta_t \frac{C_\Omega}{\eta \beta_0} \| \delta \bs{b}^n \|_0^2 
    \nonumber \\
    & & \hspace{20mm}
    \le 
    \varrho \| \delta \bs{v}_h^{0} \|_0^2
    +
    E \| \bs{D} (\delta \bs{u}_h^{0}) \|_0^2
    +
    \sum_{m=1}^n \Delta_t \frac{C_\Omega}{\eta \beta_0} \| \delta \bs{b}^m \|_0^2 
\end{eqnarray}
Dividing by $ \Delta_t^2 $ and noting $ \delta_t \bs{u}_h^n = \bs{v}_h^n $ for any $ n $, we observe that the remaining terms in Eq.~\eqref{eq_lem_stab_p2} are bounded by given data, \emph{e.g.}
\begin{equation}
    \varrho \| \delta_t \bs{v}_h^n \|_0^2      
    +
    E \| \bs{D} \bs{v}_h^n \|_0^2
    \le 
    \varrho K^2
    +
    E \| \bs{D} \bs{v}_h^{0} \|_0^2
    +
    \sum_{m=1}^n \Delta_t \frac{C_\Omega}{\eta \beta_0} \| \delta_t \bs{b}^m \|_0^2 .
\end{equation}
Combining these results, along with bounds on $ \bs{b} $ and it's derivative, we arrive at the stability estimate.
\end{pol}

\section{Optimal Prony Series Parameters for Approximating the Caputo Fractional Derivative}\label{app:timescaling}

In section \ref{meth:frac_opt}, we presented the method for optimizing the parameters of the Prony series for a defined time scale proportional to the duration of the simulation. In practice, we observed in our examples (section \ref{sect:poly_example} and \ref{sect:gao_example}) that better convergence rates can be achieved by scaling the time scale for which the Prony series parameters were mapped to. In every case, the convergence rates were significantly improved by lengthening the time scale by a factor of 10 or even 100. This was first observed as a consequence of the Geo et al. example (section \ref{sect:gao_example}), where the convergence response can be worse initially as the number of Prony terms increase. This is a results of the increase in value of $C(\bs{\beta}, \bs{\tau})$, which results in slow convergence until the approximation error hits the bounds of the error estimate in Theorem \ref{thm_prony_error}. Numerical, lengthening the time scale for mapping the Prony series parameters increases the time constants $\tau_k$ and decreases the weights $\beta_k$. We can observe from Theorem \ref{thm_prony_error} that the error bound, specific $\beta_0$ and $C(\bs{\beta}, \bs{\tau})$ reduces as a response, decreasing the overall error of the approximation. 

As this effect appears to be problem dependent (more significant in Gao et al. Example 1 \cite{gao2012finite}), we tested this in more detail using two polynomial examples: 1) the decaying oscillating function like polynomial presented in Fig. \ref{fig:polyfit}, which behaves with short variable base frequencies and 2) a polynomial that's monotonically increasing with highly exponential behavior, which behaves with single long base frequency (i.e. fitted to $e^{6*t} - 1$). The convergence response was tested for $\alpha \in \lbrace 0.1, 0.4, 0.8\rbrace$ and $N \in \lbrace 3, 6, 9, 12\rbrace$ with time scaling of 1, 10, and 100. 

The overall response is quite similar. Most important is that lengthening the time scale by a factor of 10 always results in better convergence and overall error, sometimes by two order of magnitude. Results from scaling by a factor of 100 is more variable. Most noticeably, the lower-bound with a scaling of 100 can be much worse at larger $N$. However, the convergence at larger $\Delta_t$ can be significantly better. Most noticeably, the $L^2-$norm in the example 2 can be 3 magnitude better at $\Delta_t=10^-3$ for $\alpha = 0.8$ and N = 12 with a a scaling of 100. However, also note that with a scaling of 10, the $L^2-$norm eventually surpasses the scaling of 100 after it has plateaued. Here, the gain in convergence at large $\Delta_t$ surpasses the 1 magnitude gain in the lowerbound as a result of the truncation error (Lemma 1), which requires 1.5 more magnitude in $\Delta_t$ refinement. In contrast, for example 1, we observe that the gain in convergence initially is much less significant, and the truncation error is around 2 magnitudes larger. 

Clearly, these results are also problem dependent. Part of this dependency can be observed from Theorem \ref{thm_prony_error}, by the weights of $\left [ 
    | f^\prime (0) | +
    \| f^{\prime \prime} \|_{0,1}
    \right ]$ and $ \| f \|_{W^{3,\infty}(0,T)}$ on $\epsilon$ and $C(\bs{\beta},\bs{\tau})$ respectively. 
Example 1 starts at zero, increasing quickly then oscillating until it decay exponentially to zero. Thus, $\left [ 
    | f^\prime (0) | +
    \| f^{\prime \prime} \|_{0,1}
    \right ]$ is large and $ \| f \|_{W^{3,\infty}(0,T)}$, which is the norm on the last time step is small, emphasizing the importance of $\epsilon$ in the error estimate. Thus, scaling the time scale, which improves the $C(\bs{\beta}, \bs{\tau})$ terms, sees much less benefit. In contrast, examples 2 increase exponentially from zero, resulting in $ \| f \|_{W^{3,\infty}(0,T)}$ being very large, allowing it to see more benefit from scaling the time scale. 
Another more physical reason is that example one oscillating during the duration of the simulation, with base frequencies on par with the time domain of the simulation, where as example two increase monotonically, with base frequencies much greater than the time domain of the simulation. 

In summary, the optimal Prony terms from section \ref{meth:frac_opt} can still be improved. However, the choice of further optimizing these parameters are quite subjective. The choice dependent strongly on the $\Delta_t$ used in the simulation and even behavior of the function being differentiated. We, in our examples, observed that scaling the time scaling by 10 almost always improves the results. More optimal scalings remain to be investigated. For this reason, we have chosen to a time scaling of 10 for the polynomial example (section \ref{sect:poly_example}) and 100 for the Gao \textit{et al.} example (section \ref{sect:gao_example}).


\renewcommand\thefigure{\arabic{figure}} 
\renewcommand\thetable{\arabic{table}} 
\newpage
\section*{Figures}
\vfill
\begin{figure}[h]
    \centering
    \includegraphics[width=1\textwidth]{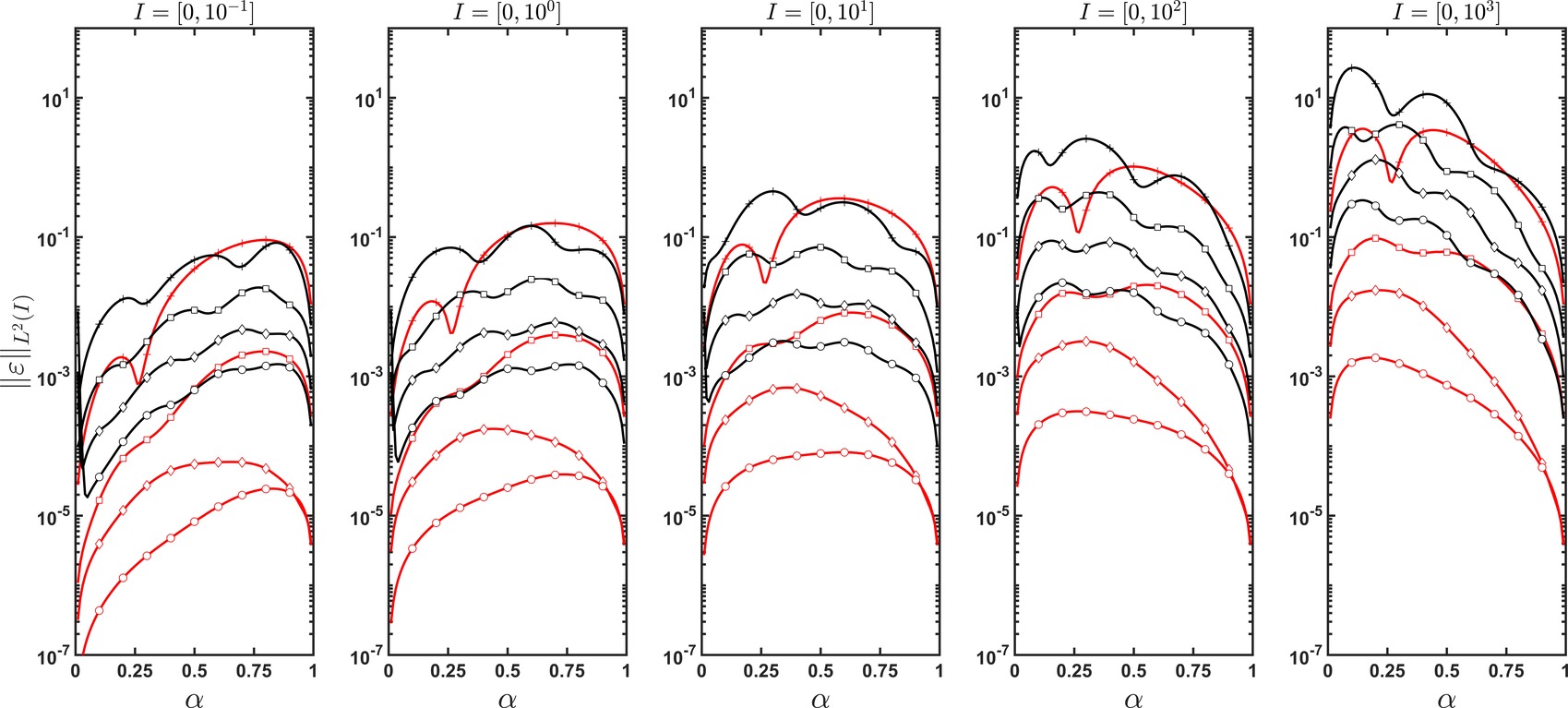} \\
    \includegraphics[width=0.35\textwidth]{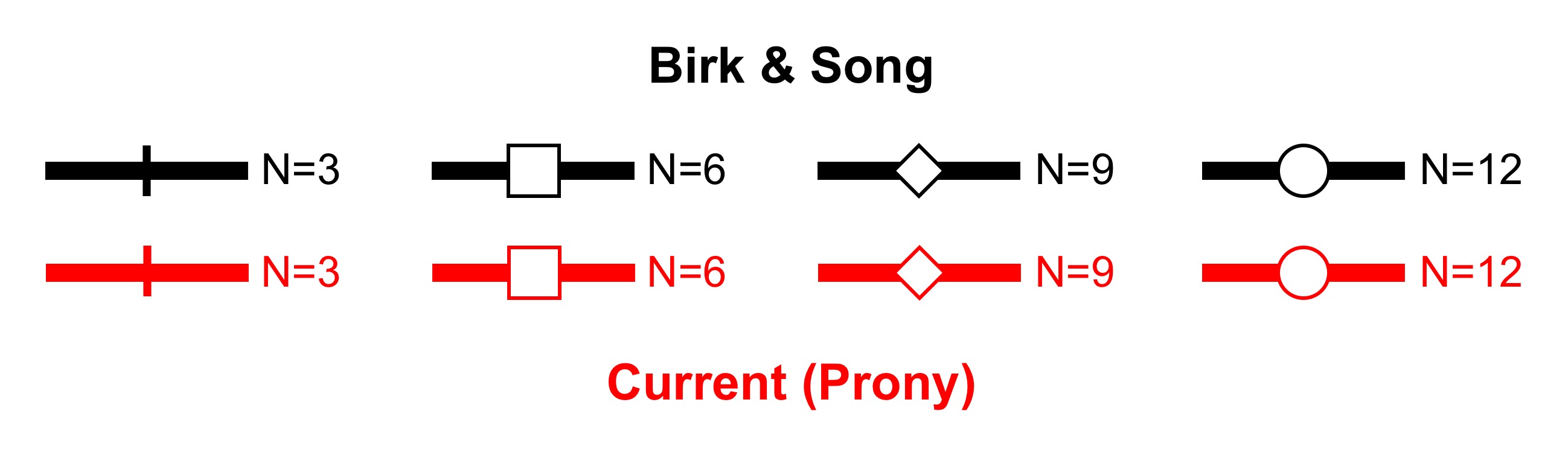} 
    \caption{Examination of the $ L^2(I)$ norm error of the truncation error $ \varepsilon $ for different time intervals, $ I $.  Each plot shows parameters selected by the optimization approach here and the integral approach in Birk and Song~\cite{birk2010improved}, for $ N \in \{3, 5, 6, 9, 12 \} $.}
    \label{fig:birkprony_err_comp}
\end{figure}
\vfill

\clearpage
\begin{figure}[h]
    \centering
    \includegraphics[width=1\textwidth]{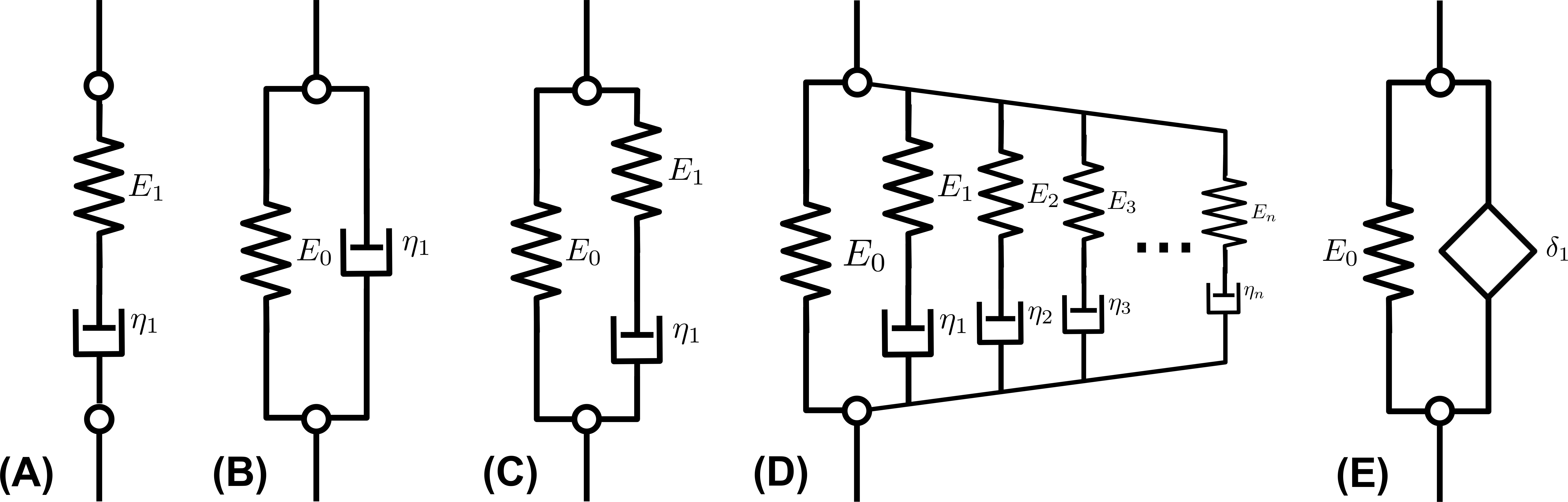} 
    \caption{Illustration of spring, dashpot and spring-pot elements used to compose viscoelastic models.  Examples include the A) Maxwell, B) Kelvin Voigt, C) Standard Linear, D) Generalized Maxwell and E) Fractional Kelvin-Voigt models.}
    \label{fig:maxandkelvinmodels}
\end{figure}

\clearpage
\begin{figure}[h]
    \centering
    \includegraphics[width=1\textwidth]{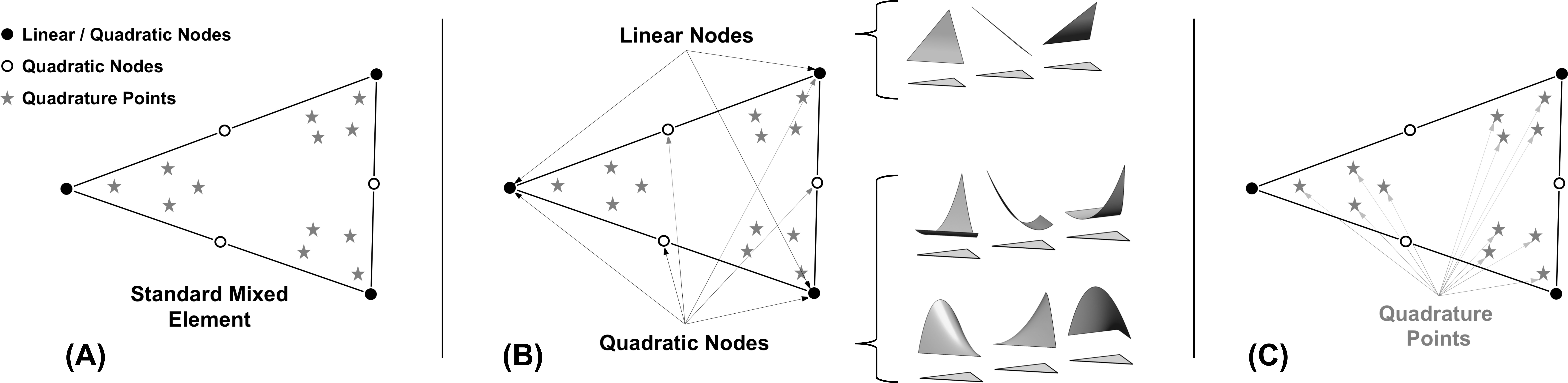} 
    \caption{Illustration of element level discretization for fractional viscoelasticity.
             (A) View of standard mixed $ \mathbb{P}^2-\mathbb{P}^1 $ triangular finite element scheme and a corresponding quadrature set.
             (B) Selection of nodes for linear and quadratic triangle along with the corresponding nodal Lagrange basis functions.
             (C) Selection of quadrature used to approximate the viscoelastic intermediate variables, $ \bs{Q} $.}
    \label{fig:fe_discretization}
\end{figure}

\clearpage
\begin{figure}
    \centering
    \includegraphics[width=\textwidth]{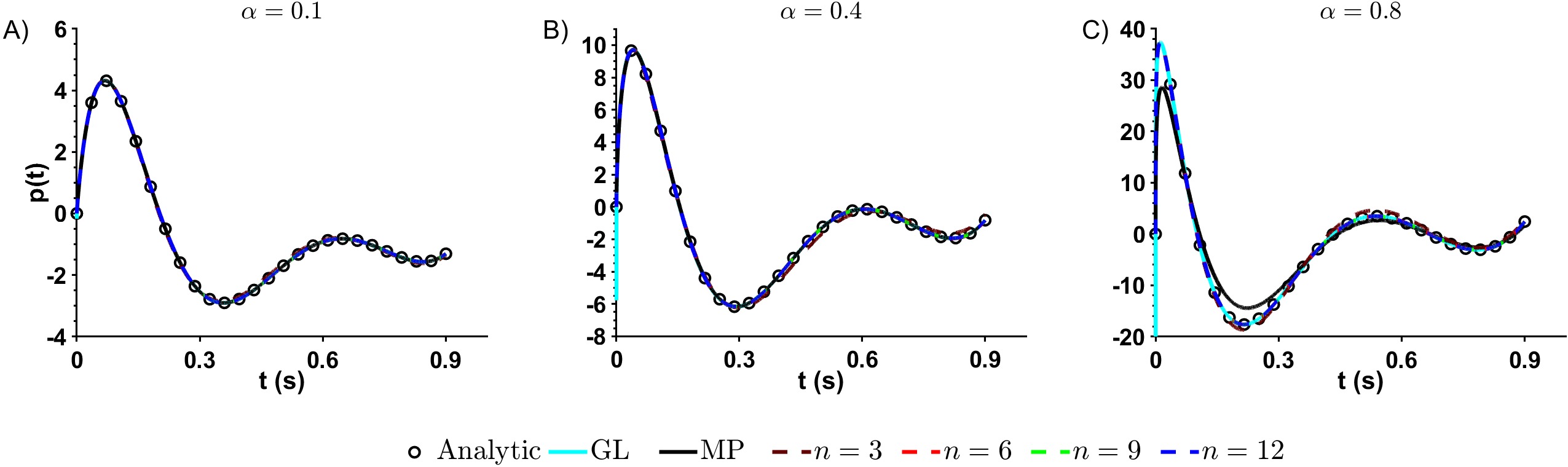}
    \caption{Examples of approximating the fractional derivative of a polynomial (with $ \{ b_1 \ldots b_8 \} = \{$2.17, 101.54, -977.47, 3368.61, -5636.44, 4937.49, -2191.59, 398.40$ \} $ ) using the midpoint rule versus the Prony-based method with $N \in \{ 3, 6, 9, 12 \} $ terms for A) $\alpha = 0.1$, B) $\alpha = 0.4$ and C) $\alpha = 0.8$.}
    \label{fig:polyfit}
\end{figure}

\clearpage
\begin{figure}
    \centering
    \includegraphics[width=\textwidth]{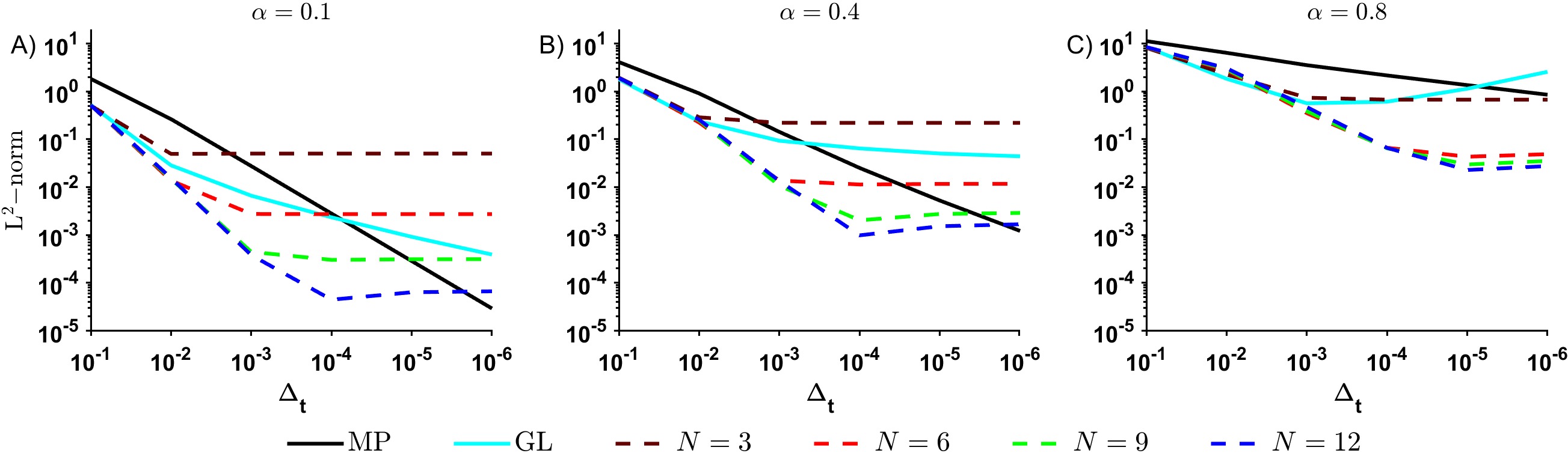} 
    \caption{The convergence of the approximation of the fractional derivative of polynomials when using the midpoint rule versus the Prony series for A) $\alpha = 0.1$, B) $\alpha = 0.4$ and C) $\alpha = 0.8$.}
    \label{fig:polyconverge}
\end{figure}

\clearpage    
\begin{figure}
    \centering
    \includegraphics[width=0.5\textwidth]{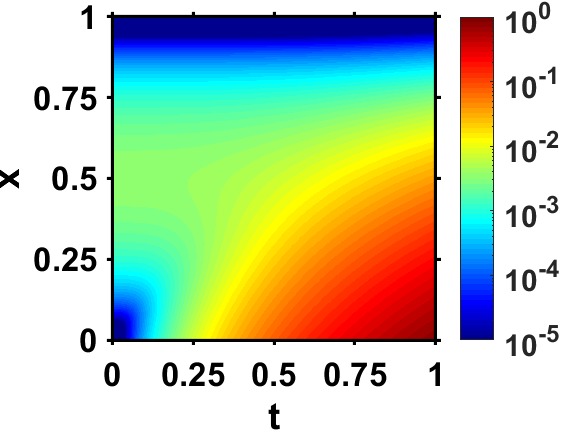}
    \caption{The solution of the Prony series approximation of the fractional differential equation example Eq. \eqref{eq:diffeqexample} with three Prony terms, $\alpha = 1/2$, $\Delta_x = 1/2500$ and $\Delta_t = 1/2500$.}
    \label{fig:FEfit}
\end{figure}

\clearpage
\begin{figure}
    \centering
    \includegraphics[width=\textwidth]{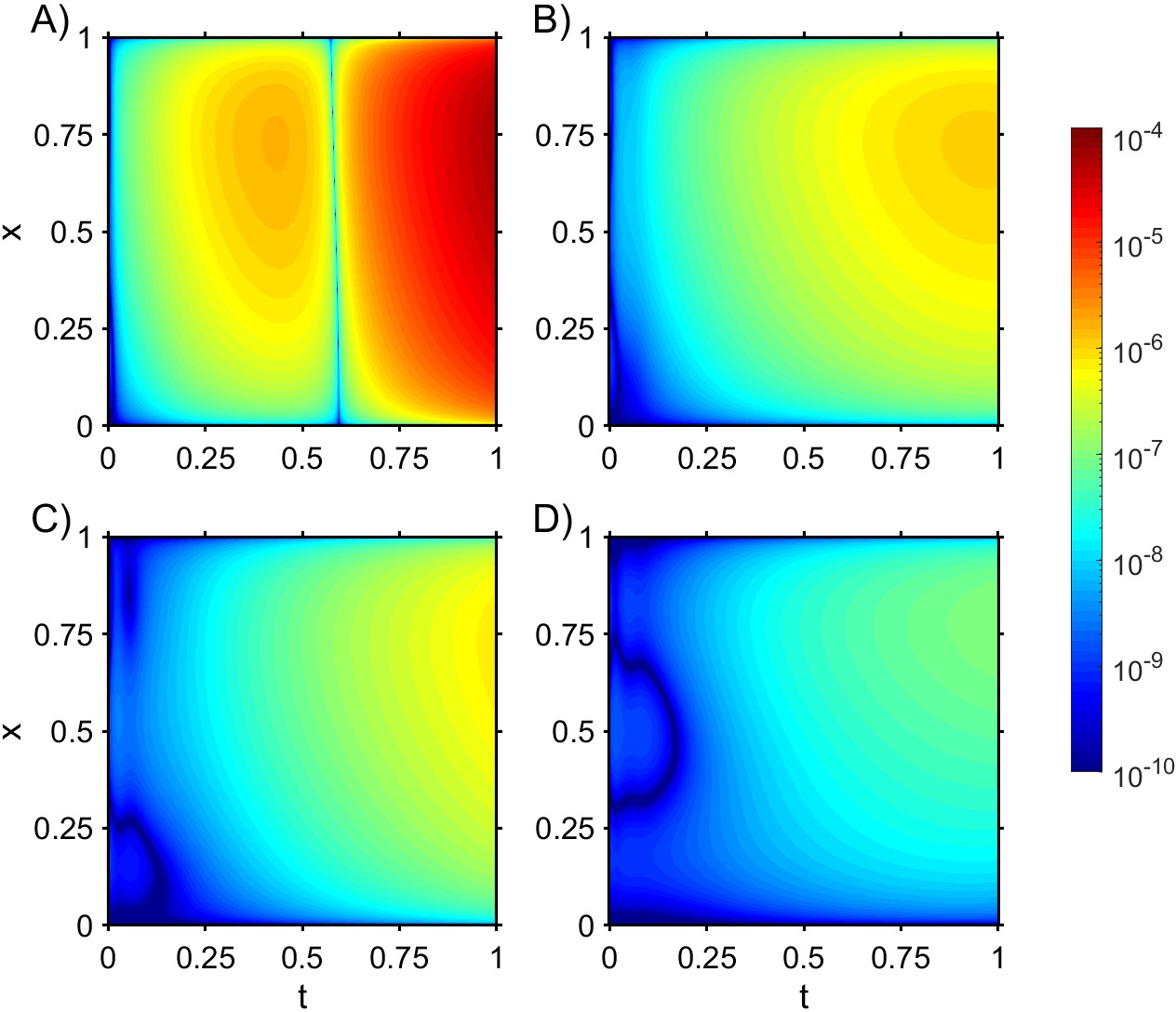}
    \caption{The absolute value of the residuals from the Prony series approximation of the fractional differential equation example Eq. \eqref{eq:diffeqexample} shown in color with $\alpha=0.5$, $\Delta_t = 1/2500$, $\Delta_x = 1/2500$, and A) 3, B) 6, C) 9, or D) 12 Prony terms in the approximation.}
    \label{fig:FEres}
\end{figure}

\clearpage
\begin{figure}
\centering
Prony-based Method, Eq.~\eqref{eq_prony_approx_final}\\
\subfloat{\includegraphics[width=\textwidth]{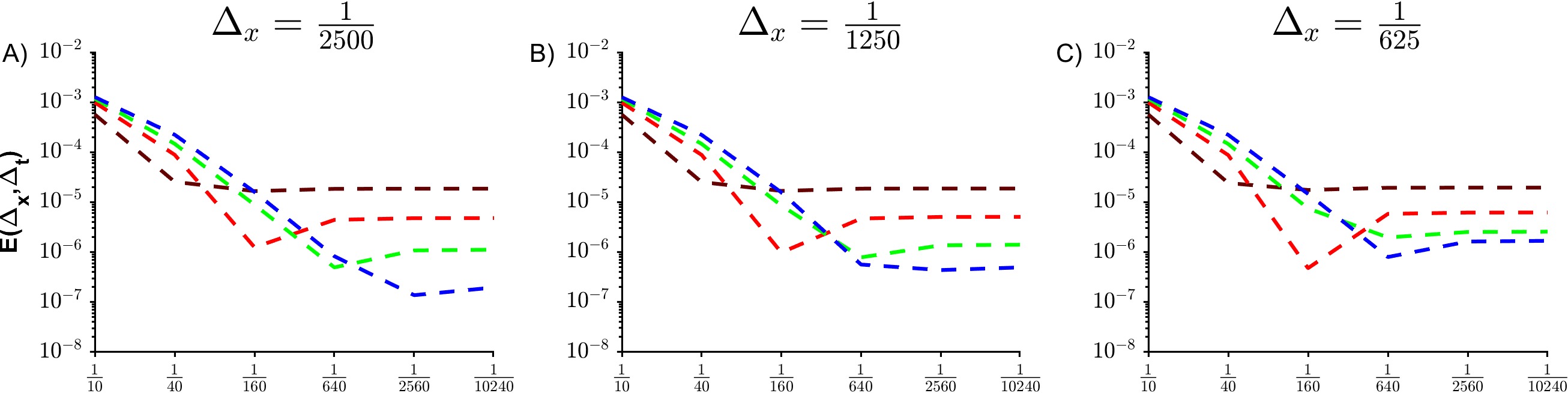}}\\
Birk and Song Method~\cite{birk2010improved}\\
\subfloat{\includegraphics[width=\textwidth]{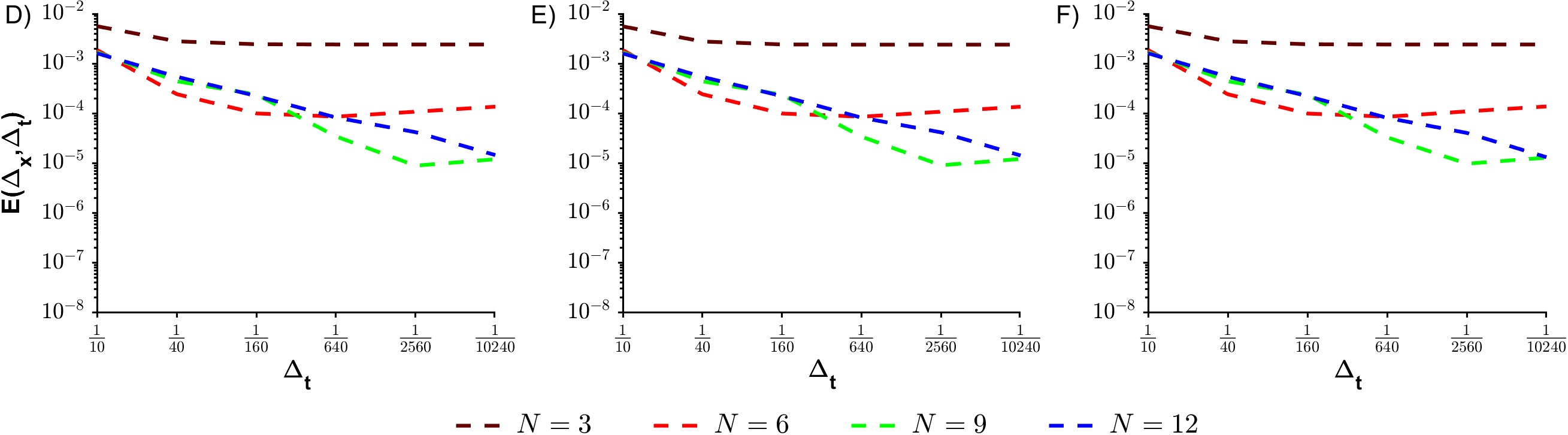}}
\caption{Convergence of the solution of the fractional diffusion equation with $\alpha=0.3$ using the Prony-based method Eq.~\eqref{eq_prony_approx_final} with refinement in time and A)$\Delta_x = 1/2500$, B) $1/1250$ and C) $1/625$ versus using the method of Birk and Song~\cite{birk2010improved} with refinement in time and D)$\Delta_x = 1/2500$, E) $1/1250$ and F) $1/625$.
}
\label{fig:gao4}
\end{figure}


\clearpage
\begin{figure}
    \centering
    \includegraphics[width=0.5\textwidth]{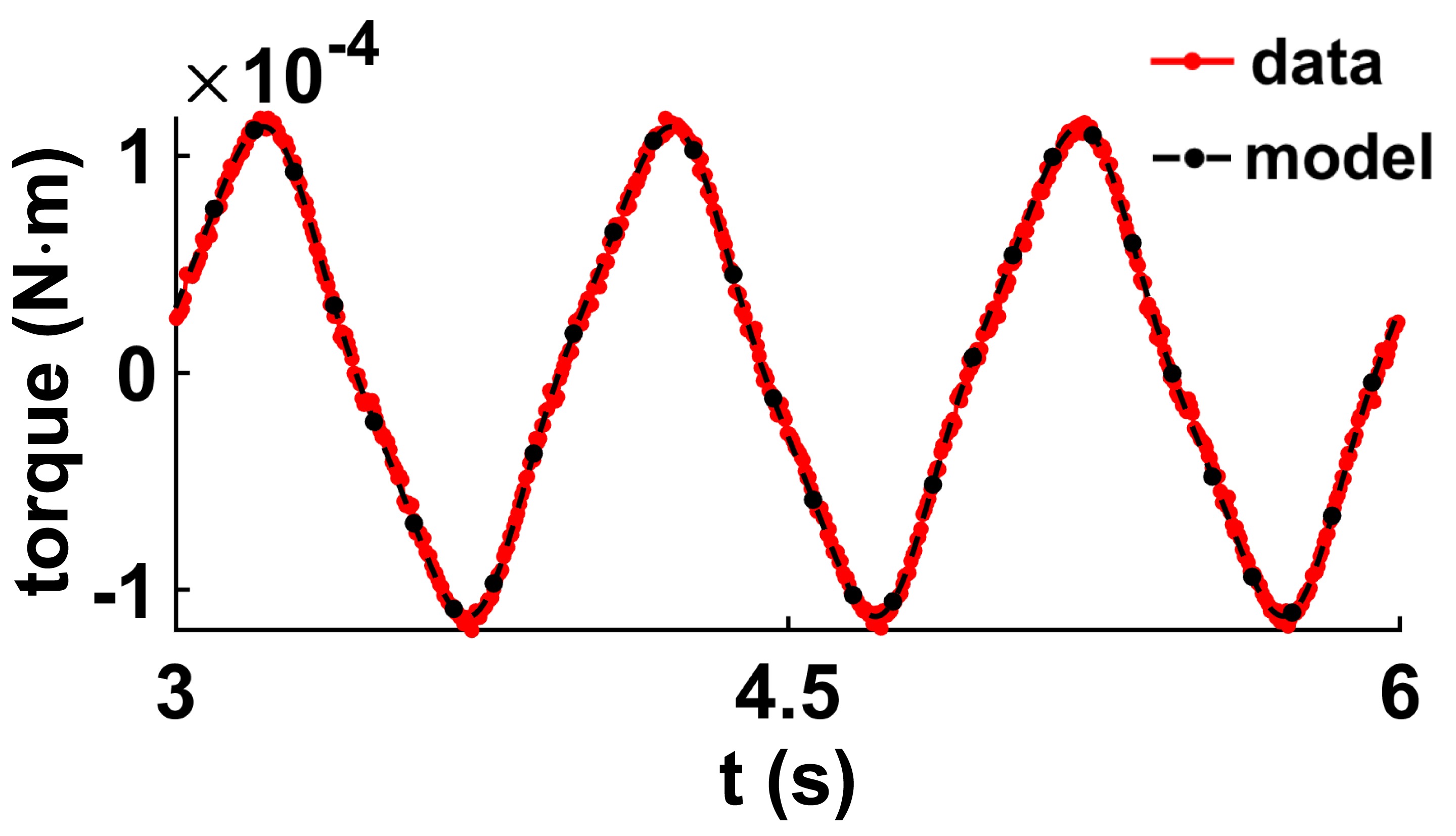}
    \caption{Comparison of the model and experimental data from Tan \emph{et al.}~\cite{Tan2013} at $10$\% compression and $25$\% shear strain, showing efficacy of the fractional viscoelastic model.}
    \label{fig:liver_behaviour}
\end{figure}

\clearpage
\begin{figure}
    \centering
    \includegraphics[width=0.9\textwidth]{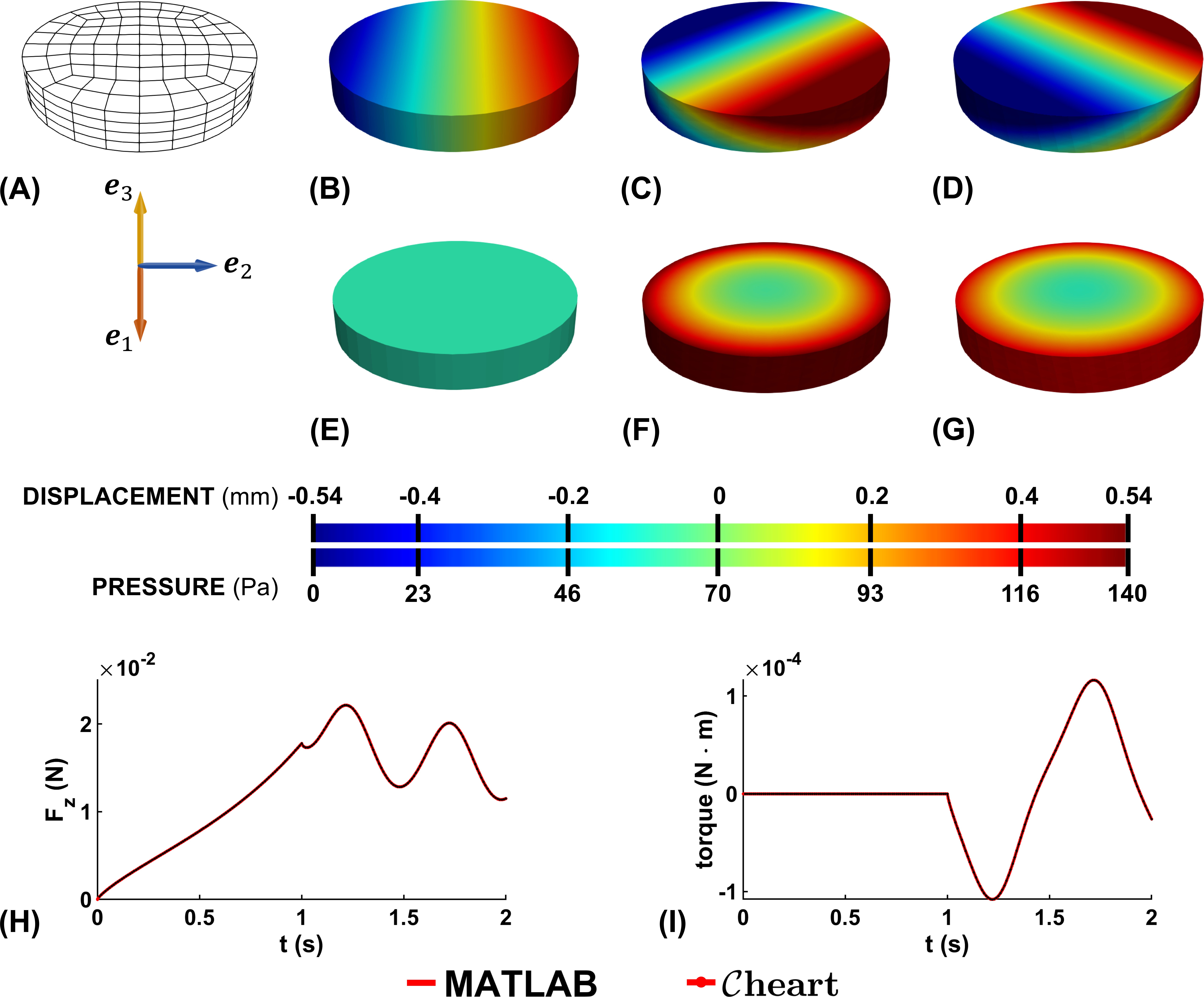}
    \caption{Illustration of the FEM solution of the liver model under deformation and shear. A) Sample mesh; $\bs{u}_2$-component of the deformation at B) compressed state ($t=1$ s) C) compressed and sheared state (counter-clockwise, $t=1.25$ s) D) compressed and sheared state (clockwise, $t=1.75$ s); Hydrostatic pressure at E) compressed state ($t=1$ s) F) compressed and sheared state (counter-clockwise, $t=1.25$ s) G) compressed and sheared state (clockwise, $t=1.75$ s); The bottom panels show the force (H) and torque response (I) of a sample under compression and shear. }
    \label{fig:liver}
\end{figure}


\clearpage
\begin{sidewaysfigure}
    \centering
    \includegraphics[width=\textwidth]{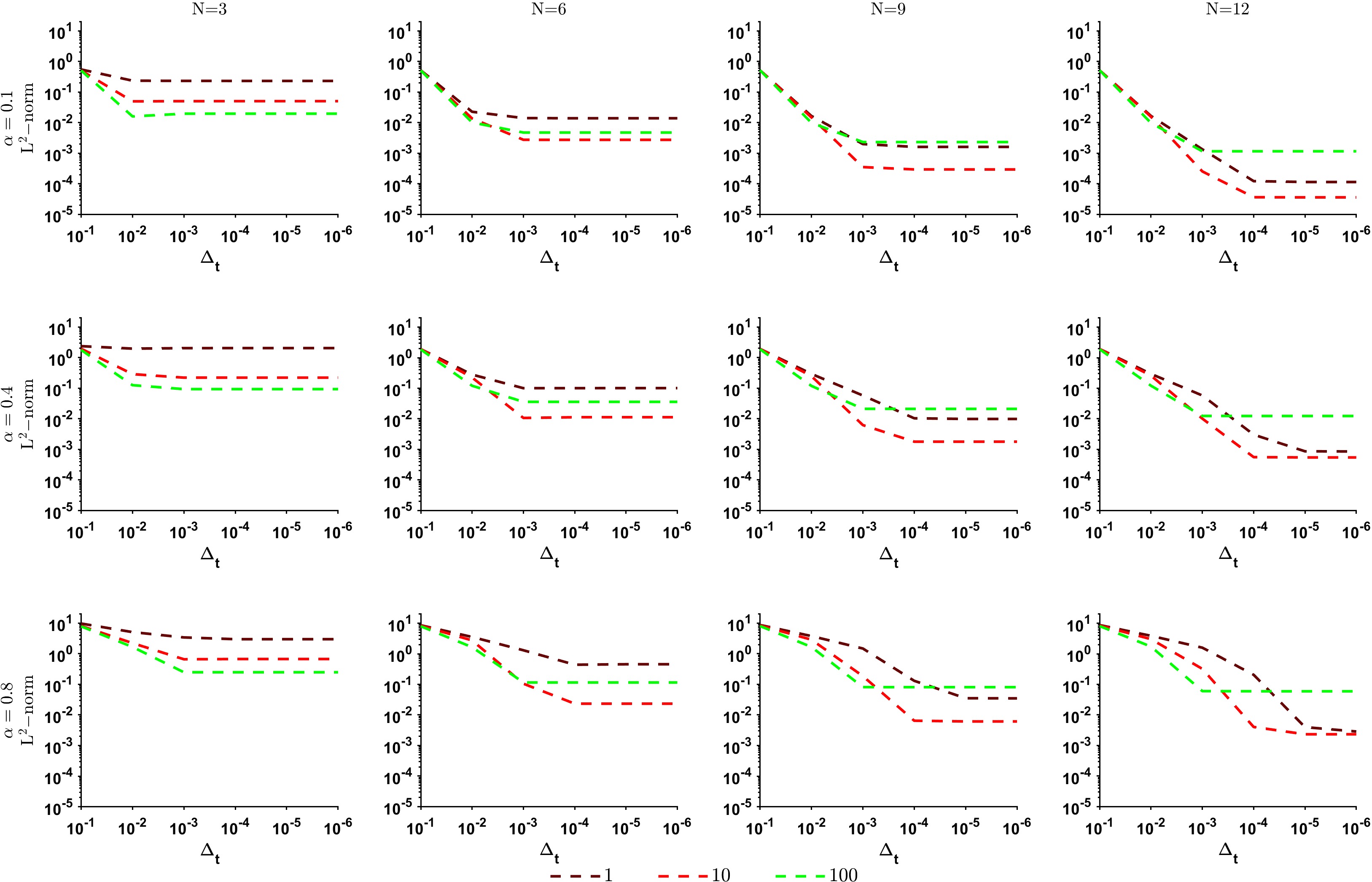}
    \caption{Comparing convergence of the Prony-based approximation Eq. \eqref{eq_prony_approx_final} in the polynomial example with refinement in time, $\alpha = 0.1$, $0.4$, $0.8$, and 3, 6, 9, and 12 Prony terms when scaling the Prony series parameters to time scales that are 1 times, 10 times and 100 times the size of the times domain.}
    \label{fig:polyscaling}
\end{sidewaysfigure}

\clearpage
\begin{sidewaysfigure}
    \centering
    \includegraphics[width=\textwidth]{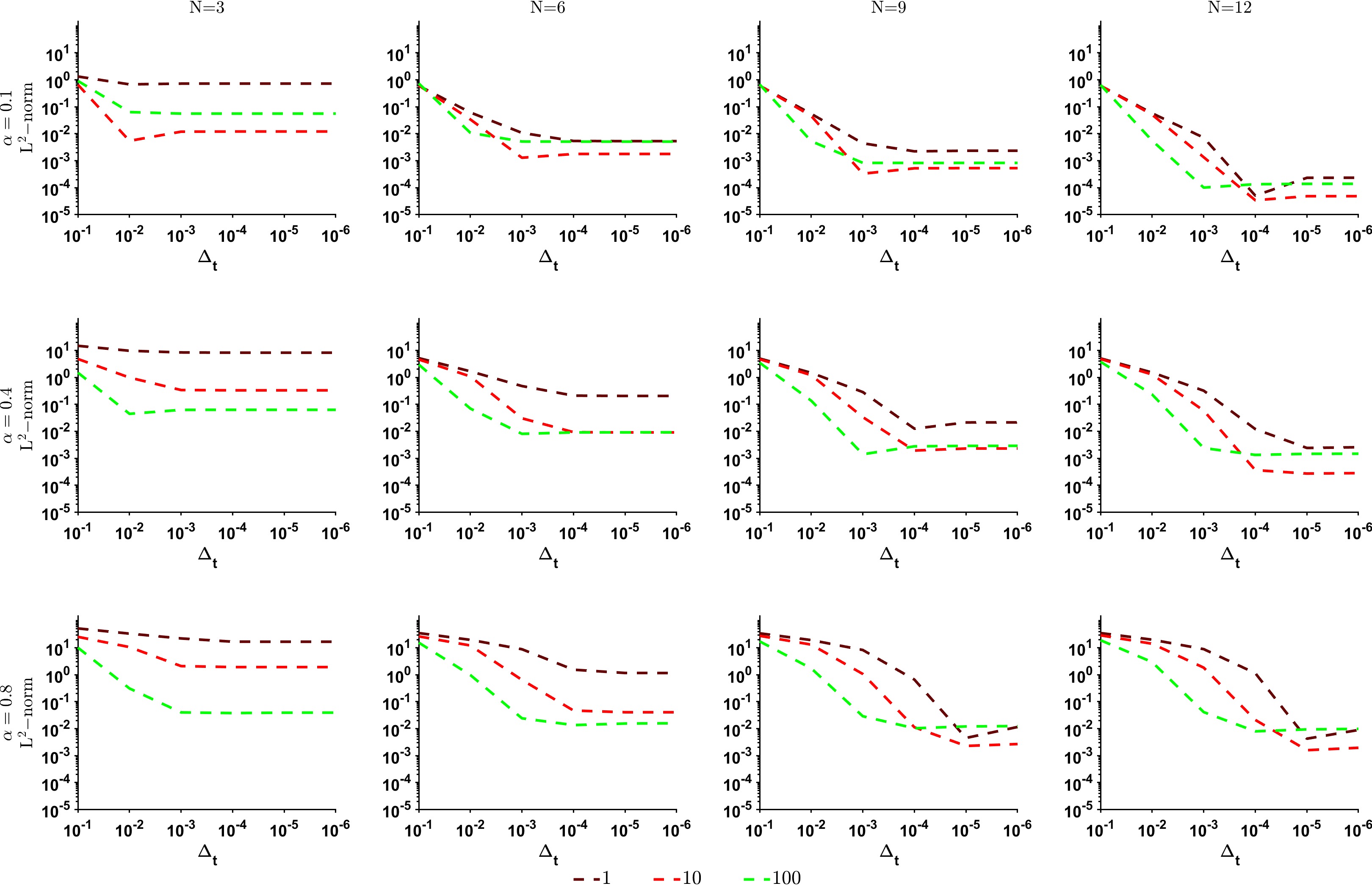}
    \caption{Comparing convergence of the Prony-based approximation Eq. \eqref{eq_prony_approx_final} in Appendix \ref{app:timescaling} exponential function example (Example 2) with refinement in time, $\alpha = 0.1$, $0.4$, $0.8$, and 3, 6, 9, and 12 Prony terms when scaling the Prony series parameters to time scales that are 1 times, 10 times and 100 times the size of the times domain.}
    \label{fig:expscaling}
\end{sidewaysfigure}

\clearpage
\section*{Tables}

\vfill

\begin{table}[h]
\centering
\begin{tabular}{R{.5in}R{.75in}R{.5in}R{.5in}R{.5in}R{.5in}}
\hline
 \multicolumn{1}{c}{\textbf{$\Delta_t$} (s)} & \multicolumn{1}{c}{\textbf{Midpoint} (s)} & \multicolumn{4}{c}{\textbf{Prony} (s)}   \\
\hline
 \multicolumn{1}{c}{}  & \multicolumn{1}{c}{} & \multicolumn{1}{c}{$N = 3$} & \multicolumn{1}{c}{$N = 6$} & \multicolumn{1}{c}{$N = 9$} & \multicolumn{1}{c}{$N = 12$} \\
\hline \hline
1e-4	&8.16	&0.28	&0.27	&0.30	&0.33\\
5e-5	&29.25	&0.53	&0.58	&0.60	&0.61\\
1e-5	&785.14	&2.58	&2.64	&2.87	&3.01\\
\hline
\end{tabular}
\caption{Computational cost of approximating the fractional derivative of the polynomial in Eq.~\eqref{eqn:polynomialfunction} with the midpoint rule versus the Prony-based method with $ N \in \{ 3, 6, 9, 12 \} $.  Here the time interval for the calculation is $ [0,0.9] $ and $ \alpha = 0.4$.}\label{tab:polynomialtime}
\end{table}

\vfill

\clearpage
\begin{table}
\setlength\tabcolsep{3pt}
\centering
\begin{tabular}{C{.1in} R{.3in}  C{.01in} L{.45in} C{.01in} L{.45in} C{.01in} L{.45in} L{.45in} L{.45in} L{.45in} C{.01in} L{.45in} L{.45in} L{.45in} L{.45in}}
\multicolumn{1}{c}{} & \multicolumn{1}{c}{\textbf{$\Delta_t$}} & \multicolumn{1}{c}{} & \multicolumn{1}{c}{\textbf{MP}} & \multicolumn{1}{c}{} & \multicolumn{1}{c}{\textbf{GL}} & \multicolumn{1}{c}{} & \multicolumn{4}{c}{\textbf{Prony}} & \multicolumn{1}{c}{} & \multicolumn{4}{c}{\textbf{Birk and Song}} \\
\cmidrule{1-2} \cmidrule{4-4} \cmidrule{6-6} \cmidrule{8-11} \cmidrule{13-16}
\multicolumn{1}{c}{} & \multicolumn{1}{c}{}  & \multicolumn{1}{c}{}  & \multicolumn{1}{c}{} & \multicolumn{1}{c}{}  & \multicolumn{1}{c}{}  & \multicolumn{1}{c}{} &\multicolumn{1}{c}{$N = 3$} & \multicolumn{1}{c}{$N = 6$} & \multicolumn{1}{c}{$N = 9$} & \multicolumn{1}{c}{$N = 12$}   & \multicolumn{1}{c}{} & \multicolumn{1}{c}{$N = 3$} & \multicolumn{1}{c}{$N = 6$} & \multicolumn{1}{c}{$N = 9$} & \multicolumn{1}{c}{$N = 12$} \\
\cmidrule{1-2} \cmidrule{4-4} \cmidrule{6-6} \cmidrule{8-11} \cmidrule{13-16}
\morecmidrules
\cmidrule{1-2} \cmidrule{4-4} \cmidrule{6-6} \cmidrule{8-11} \cmidrule{13-16}
\parbox[t]{0.1in}{\multirow{6}{*}{\rotatebox[origin=c]{90}{$\alpha=0.1$}}} 
&1e-1  &  &1.84e0  &  &5.24e-1  &  &5.07e-1  &5.14e-1  &5.14e-1  &5.14e-1  &  &4.73e-1  &4.95e-1  &5.12e-1  &5.14e-1\\
&1e-2  &  &2.64e-1  &  &2.86e-2  &  &4.96e-2  &1.37e-2  &1.49e-2  &1.56e-2  &  &2.99e-1  &6.78e-2  &2.77e-2  &1.84e-2\\
&1e-3  &  &2.76e-2  &  &6.64e-3  &  &5.04e-2  &2.76e-3  &4.45e-4  &3.81e-4  &  &3.01e-1  &7.42e-2  &2.23e-2  &7.31e-3\\
&1e-4  &  &2.81e-3  &  &2.35e-3  &  &5.04e-2  &2.74e-3  &3.03e-4  &4.47e-5  &  &3.01e-1  &7.46e-2  &2.24e-2  &6.81e-3\\
&1e-5  &  &2.86e-4  &  &9.19e-4  &  &5.04e-2  &2.75e-3  &3.13e-4  &6.39e-5  &  &3.01e-1  &7.46e-2  &2.24e-2  &6.78e-3\\
&1e-6  &  &2.93e-5  &  &3.92e-4  &  &5.04e-2  &2.75e-3  &3.14e-4  &6.69e-5  &  &3.01e-1  &7.46e-2  &2.24e-2  &6.77e-3\\
\cmidrule{1-2} \cmidrule{4-4} \cmidrule{6-6} \cmidrule{8-11} \cmidrule{13-16}
&0     &  &- -      &  &- -      &  &5.04e-2  &2.75e-3  & 3.14e-4  & 6.71e-5  &  &3.01e-1  &7.46e-2  &2.24e-2  &6.77e-3\\ 
\cmidrule{1-2} \cmidrule{4-4} \cmidrule{6-6} \cmidrule{8-11} \cmidrule{13-16}
\morecmidrules
\cmidrule{1-2} \cmidrule{4-4} \cmidrule{6-6} \cmidrule{8-11} \cmidrule{13-16}
\parbox[t]{0.1in}{\multirow{6}{*}{\rotatebox[origin=c]{90}{$\alpha=0.4$}}} 
&1e-1  &  &4.12e0  &  &1.80e0  &  &1.95e0  &1.91e0  &1.91e0  &1.92e0  &  &2.73e0  &2.22e0  &1.99e0  &1.94e0\\
&1e-2  &  &9.12e-1  &  &2.36e-1  &  &2.88e-1  &2.25e-1  &2.48e-1  &2.58e-1  &  &1.23e0  &5.25e-1  &3.34e-1  &3.22e-1\\
&1e-3  &  &1.43e-1  &  &9.34e-2  &  &2.21e-1  &1.38e-2  &1.08e-2  &1.35e-2  &  &1.21e0  &4.74e-1  &1.63e-1  &7.59e-2\\
&1e-4  &  &2.53e-2  &  &6.46e-2  &  &2.21e-1  &1.14e-2  &2.04e-3  &9.88e-4  &  &1.21e0  &4.17e-1  &1.43e-1  &4.39e-2\\
&1e-5  &  &5.27e-3  &  &5.05e-2  &  &2.21e-1  &1.17e-2  &2.76e-3  &1.53e-3  &  &1.21e0  &4.12e-1  &1.38e-1  &4.31e-2\\
&1e-6  &  &1.22e-3  &  &4.42e-2  &  &2.21e-1  &1.18e-2  &2.90e-3  &1.68e-3  &  &1.21e0  &4.12e-1  &1.37e-1  &4.34e-2\\
\cmidrule{1-2} \cmidrule{4-4} \cmidrule{6-6} \cmidrule{8-11} \cmidrule{13-16}
&0     &  & - -      &  & - -      &  &2.21e-1  & 1.18e-2  &2.92e-3  &1.69e-3  &  &1.21e0  &4.12e-1  &1.37e-1  &4.35e-2 \\ 
\cmidrule{1-2} \cmidrule{4-4} \cmidrule{6-6} \cmidrule{8-11} \cmidrule{13-16}
\morecmidrules
\cmidrule{1-2} \cmidrule{4-4} \cmidrule{6-6} \cmidrule{8-11} \cmidrule{13-16}
\parbox[t]{0.1in}{\multirow{6}{*}{\rotatebox[origin=c]{90}{$\alpha=0.8$}}} 
&1e-1  &  &1.13e1  &  &8.46e0  &  &8.14e0  &8.37e0  &8.42e0  &8.47e0  &  &8.22e0  &9.10e0  &9.38e0  &9.36e0\\
&1e-2  &  &6.42e0  &  &1.84e0  &  &2.27e0  &2.79e0  &2.92e0  &3.06e0  &  &4.54e0  &5.35e0  &5.90e0  &5.53e0\\
&1e-3  &  &3.54e0  &  &5.66e-1  &  &7.41e-1  &3.52e-1  &3.86e-1  &4.71e-1  &  &4.48e0  &2.80e0  &3.19e0  &3.70e0\\
&1e-4  &  &2.16e0  &  &6.11e-1  &  &6.75e-1  &6.67e-2  &6.50e-2  &6.62e-2  &  &3.98e0  &1.71e0  &2.41e0  &2.07e0\\
&1e-5  &  &1.36e0  &  &1.14e0  &  &6.75e-1  &4.34e-2  &2.98e-2  &2.29e-2  &  &2.85e0  &1.70e0  &1.13e0  &1.42e0\\
&1e-6  &  &8.56e-1  &  &2.58e0  &  &6.75e-1  &4.90e-2  &3.54e-2  &2.76e-2  &  &2.89e0  &1.25e0  &1.10e0  &7.39e-1\\
\cmidrule{1-2} \cmidrule{4-4} \cmidrule{6-6} \cmidrule{8-11} \cmidrule{13-16}
&0     &  &- -       &  & - -      &  & 6.75e-1  & 4.98e-2  & 3.64e-2  & 2.87e-2  &  &2.89e0  &8.03e-1 &2.63e-1 &8.89e-2 \\ 
\cmidrule{1-2} \cmidrule{4-4} \cmidrule{6-6} \cmidrule{8-11} \cmidrule{13-16}
\morecmidrules
\cmidrule{1-2} \cmidrule{4-4} \cmidrule{6-6} \cmidrule{8-11} \cmidrule{13-16}
\end{tabular}
\caption{Comparing the approximation error for the fractional derivative of polynomials when using the midpoint (MP) rule versus the Prony series.
Here $ \Delta_t = 0$ denotes the norm, $\| \varepsilon \|_0$, of the truncation error predicted in Lemma~\ref{lem_Prony_error}. }\label{tab:polynomialconverge}
\end{table}

\clearpage
\begin{table}
\small
\setlength\tabcolsep{2pt}
\centering
\begin{tabular}{C{.2in}L{.15in}L{.35in}R{.5in}R{.28in}C{0.01in}R{.5in}R{.28in}R{.5in}R{.28in}R{.5in}R{.28in}R{.5in}R{.28in}}
\cline{4-5} \cline{7-14}
\multicolumn{3}{c}{}& \multicolumn{2}{c}{\textbf{Gao et al.} \cite{gao2012finite}} & & \multicolumn{8}{c}{\textbf{Prony approximation}} \\
\cline{4-5} \cline{7-14}
\multicolumn{5}{c}{} & & \multicolumn{2}{c}{$N$ = 3} & \multicolumn{2}{c}{$N$ = 6} & \multicolumn{2}{c}{$N$ = 9} & \multicolumn{2}{c}{$N$ = 12}\\
\cline{1-5} \cline{7-14}
$\alpha$& \multicolumn{1}{c}{$h$} & \multicolumn{1}{c}{$\Delta_t$} & \multicolumn{1}{c}{\scriptsize E($h$,$\Delta_t$)} & log$_2$ & & \multicolumn{1}{c}{\scriptsize E($h$,$\Delta_t$)} & log$_2$ & \multicolumn{1}{c}{\scriptsize E($h$,$\Delta_t$)} & log$_2$ & \multicolumn{1}{c}{\scriptsize E($h$,$\Delta_t$)} & log$_2$ & \multicolumn{1}{c}{\scriptsize E($h$,$\Delta_t$)} & log$_2$\\
\cline{1-5} \cline{7-14}
\parbox[t]{0.2in}{\multirow{6}{*}{1/2}} &
\parbox[t]{0.15in}{\multirow{6}{*}{\rotatebox[origin=c]{90}{$1/20\,000$}}}
   &1/10 	&8.48e-4	&1.41	& &2.13e-3	&1.62	&3.33e-3	&1.26	&3.73e-3	&0.96	&3.96e-3	&0.74\\
   & &1/20 	&3.19e-4	&1.44	& &6.94e-4	&1.54	&1.39e-3	&1.73	&1.91e-3	&1.53	&2.36e-3	&1.25\\
   & &1/40 	&1.17e-4	&1.46	& &2.39e-4	&1.03	&4.17e-4	&1.94	&6.64e-4	&1.86	&9.94e-4	&1.73\\
   & &1/80 	&4.26e-5	&1.47	& &1.17e-4	&0.44	&1.09e-4	&2.03	&1.83e-4	&1.98	&3.00e-4	&1.92\\
   & &1/160	&1.53e-5	&1.48	& &8.60e-5	&0.14	&2.67e-5	&2.21	&4.64e-5	&2.08	&7.92e-5	&1.99\\
   & &1/320	&5.49e-6	&*	& &7.82e-5	&*	&5.79e-6	&*	&1.10e-5	&*	&2.00e-5	&*\\
\cline{1-5} \cline{7-14}
\parbox[t]{0.2in}{\multirow{6}{*}{2/3}} &
\parbox[t]{0.15in}{\multirow{6}{*}{\rotatebox[origin=c]{90}{$1/20\,000$}}}
   & 1/10 	& 1.91e-3	&1.26	& & 4.52e-3	& 1.52	& 6.78e-3	& 0.99	& 7.40e-3	& 0.72	& 7.88e-3	& 0.56\\
   & & 1/20 	& 7.96e-4	& 1.29	& & 1.58e-3	& 1.65	& 3.42e-3	& 1.58	& 4.51e-3	& 1.31	& 5.33e-3	& 1.00\\
   & & 1/40 	& 3.25e-4	& 1.31	& & 5.04e-4	& 1.31	& 1.14e-3	& 1.90	& 1.82e-3	& 1.78	& 2.67e-3	& 1.59\\
   & & 1/80 	& 1.31e-4	& 1.32	& & 2.04e-4	& 0.69	& 3.07e-4	& 2.07	& 5.28e-4	& 2.01	& 8.85e-4	& 1.91\\
   & & 1/160	& 5.26e-5	& 1.32	& & 1.26e-4	& 0.24	& 7.28e-5	& 2.53	& 1.31e-4	& 2.32	& 2.35e-4	& 2.13\\
   & & 1/320	& 2.10e-5	& *	& & 1.07e-4	& *	& 1.26e-5	& *	& 2.64e-5	& *	& 5.36e-5	& *\\
\cline{1-5} \cline{7-14}
\end{tabular}
\caption{Reproducing Table 1 from Gao \textit{et al.} \cite{gao2012finite} using the Prony-based method Eq. \eqref{eq_prony_approx_final} for refinement in time with a constant $h = $ 5e-5.}\label{tab:gaotable1}
\end{table}

\clearpage
\begin{table}
\small
\setlength\tabcolsep{2pt}
\centering
\begin{tabular}{C{.2in}L{.35in}L{.15in}R{.5in}R{.28in}C{0.01in}R{.5in}R{.28in}R{.5in}R{.28in}R{.5in}R{.28in}R{.5in}R{.28in}}
\cline{4-5} \cline{7-14}
\multicolumn{3}{c}{}& \multicolumn{2}{c}{\textbf{Gao et al.} \cite{gao2012finite}} & & \multicolumn{8}{c}{\textbf{Prony approximation}} \\
\cline{4-5} \cline{7-14}
\multicolumn{5}{c}{} & & \multicolumn{2}{c}{$N$ = 3} & \multicolumn{2}{c}{$N$ = 6} & \multicolumn{2}{c}{$N$ = 9} & \multicolumn{2}{c}{$N$ = 12}\\
\cline{1-5} \cline{7-14}
$\alpha$& \multicolumn{1}{c}{$h$} &\multicolumn{1}{c}{$\Delta_t$} &  \multicolumn{1}{c}{\scriptsize E($h$,$\Delta_t$)} & log$_2$ & & \multicolumn{1}{c}{\scriptsize E($h$,$\Delta_t$)} & log$_2$ & \multicolumn{1}{c}{\scriptsize E($h$,$\Delta_t$)} & log$_2$ & \multicolumn{1}{c}{\scriptsize E($h$,$\Delta_t$)} & log$_2$ & \multicolumn{1}{c}{\scriptsize E($h$,$\Delta_t$)} & log$_2$\\
\cline{1-5} \cline{7-14}
\parbox[t]{0.2in}{\multirow{5}{*}{1/2}}
   &1/10 &\parbox[t]{0.15in}{\multirow{5}{*}{\rotatebox[origin=c]{90}{$1/20\,000$}}} 	&1.50e-2	&2.00	& &5.95e-3	&2.06	&6.02e-3	&2.01	&6.02e-3	&2.01	&6.01e-3	&2.01\\
   &1/20 	& &3.75e-3	&2.00	& &1.43e-3	&2.22	&1.50e-3	&2.00	&1.50e-3	&2.00	&1.50e-3	&2.00\\
   &1/40 	& &9.38e-4  &2.00	& &3.07e-4	&3.24	&3.75e-4	&1.99	&3.75e-4	&1.99	&3.74e-4	&2.00\\
   &1/80 	& &2.35e-4	&2.00	& &3.25e-5	&-0.75	&9.47e-5	&1.95	&9.44e-5	&1.96	&9.37e-5	&2.00\\
   &1/160	& &5.86e-5	&*	& &5.47e-5	&*	&2.45e-5	&*	&2.43e-5	&*	&2.35e-5	&*\\
\cline{1-5} \cline{7-14}
\parbox[t]{0.2in}{\multirow{5}{*}{2/3}}
   &1/10 &\parbox[t]{0.15in}{\multirow{5}{*}{\rotatebox[origin=c]{90}{$1/20\,000$}}} 	&1.37e-2	&1.99	& &5.67e-3	&2.08	&5.77e-3	&2.00	&5.77e-3	&2.00	&5.77e-3	&2.00\\
   &1/20 	& &3.44e-3	&2.00	& &1.35e-3	&2.31	&1.44e-3	&1.98	&1.44e-3	&1.97	&1.44e-3	&1.98\\
   &1/40 	& &8.61e-4	&2.00	& &2.71e-4	&3.33	&3.66e-4	&1.92	&3.67e-4	&1.90	&3.67e-4	&1.91\\
   &1/80 	& &2.15e-4	&2.00	& &2.68e-5	&-1.58	&9.67e-5	&1.71	&9.81e-5	&1.66	&9.75e-5	&1.68\\
   &1/160	& &5.37e-5	&*	& &8.01e-5	&*	&2.96e-5	&*	&3.11e-5	&*	&3.05e-5	&*\\
\cline{1-5} \cline{7-14}
\end{tabular}
\caption{Reproducing Table 2 from Gao \textit{et al.} \cite{gao2012finite} using the Prony-based method Eq. \eqref{eq_prony_approx_final} for refinement in space with a constant $\Delta_t =$ 5e-5.}\label{tab:gaotable2}
\end{table}

\clearpage
\begin{table}
\centering
\begin{tabular}{R{.25in}|R{.5in}|R{.5in}|R{.5in}|R{.5in}|R{.5in}|R{.5in}}
 \hline
 \multicolumn{1}{c}{}  &  \multicolumn{1}{c}{\textbf{$\Delta_t$} (s)} &  \multicolumn{1}{c}{\textbf{GL} (s)} & \multicolumn{4}{c}{\textbf{Prony} (s)}   \\
 \hline 
 \multicolumn{1}{c}{} & \multicolumn{1}{c}{}  & \multicolumn{1}{c}{} & \multicolumn{1}{c}{$N = 3$} & \multicolumn{1}{c}{$N = 6$} & \multicolumn{1}{c}{$N = 9$} & \multicolumn{1}{c}{$N = 12$} \\
\hline \hline
\parbox[t]{0.15in}{\multirow{7}{*}{\rotatebox[origin=c]{90}{$\alpha=0.2$ }}} 
	&1e-2	 & 0.029	& 0.054  & 0.044 & 0.051 & 0.047 \\
	&5e-3	 & 0.052	& 0.086	 & 0.093 & 0.097 & 0.094 \\
	&1e-3	 & 0.202	& 0.406	 & 0.429 & 0.442 & 0.561 \\
	&5e-4    & 0.492	& 0.818	 & 0.871 & 0.890 & 0.933 \\
	&1e-4    & 5.728	& 4.676	 & 4.169 & 4.348 & 4.530 \\
	&5e-5    & 17.706	& 8.392	 & 8.369 & 8.859 & 9.088 \\
	&1e-5    & 770.569	& 38.051 & 37.888 & 39.370 & 41.773 \\
\hline
\end{tabular}
\caption{Comparison of solving speed for the Gr{\"{u}}nwald-Letnikov (GL) and Prony approximation methods (MATLAB implementation). 
The time results reported correspond to calculations for a single point on the outer edge of the cylinder.
All methods were run in serial on a DELL PRECISION M4800 quad-core Intel(R) Core(TM) i7-4810MQ CPU @ 2.80 GHz.}\label{tab:oldnewtimecomparison}
\end{table}

\clearpage
\begin{table}
\centering
\begin{tabular}{R{.5in}|R{.5in}|R{.5in}|R{.5in}|R{.5in}}
\hline
  \multicolumn{1}{c}{\textbf{$\Delta_t$} (s)} &  \multicolumn{4}{c}{\textbf{Prony error} (\%)}   \\
 \hline 
 \multicolumn{1}{c}{} & \multicolumn{1}{c}{$N = 3$} & \multicolumn{1}{c}{$N = 6$} & \multicolumn{1}{c}{$N = 9$} & \multicolumn{1}{c}{$N = 12$} \\
\hline \hline
1e-2	 	& 2.26e0    & 5.05e-1 & 6.51e-1 & 7.01e-1 \\
5e-3		& 2.06e0	& 1.52e-1 & 2.50e-1 & 3.27e-1 \\
1e-3	 	& 2.00e0	& 1.01e-1 & 1.30e-2 & 1.97e-2 \\
5e-4    	& 1.99e0	& 1.04e-1 & 8.79e-3 & 4.86e-3 \\
\hline
\end{tabular}
\caption{Accuracy and convergence of the Prony approximation with $\Delta_t$ and $N$ applied to the ideal analytic solution.
Percentage errors were measured relative to the GL with a refined time-step size of $\Delta_t=$ 1e-5.
The percentage errors reported in the table correspond to the $L^2(0,T)$ norm difference in (1,3)-Cauchy stress component normalized by the $L^2(0,T)$ norm of the (1,3)-Cauchy stress component predicted by the reference solution.
}\label{tab:newaccuracy}
\end{table}

\clearpage
\begin{table}
\centering
\begin{tabular}{R{1.0in}|R{.5in}|R{.5in}|R{.5in}}
 \hline
 \textbf{Model} $\backslash \Delta_t$ & 1e-2 & 1e-3 & 1e-4 \\
 \hline \hline
 \textbf{Hyperelastic} & 170.43s & 672.71s & 2899.40s \\
 \textbf{Viscoelastic} & 190.24s & 669.65s & 3365.30s \\
\hline
\end{tabular}
\hspace{10mm}
\begin{tabular}{R{.5in}|R{.5in}|R{.5in}|R{.5in}}
\hline
 \multicolumn{4}{c}{\textbf{Viscoelastic}}    \\
 \hline 
 \multicolumn{1}{c}{$N = 3$} & \multicolumn{1}{c}{$N = 6$} & \multicolumn{1}{c}{$N = 9$} & \multicolumn{1}{c}{$N = 12$} \\
 \hline \hline
  655.69s & 684.74s & 669.65s & 677.85s \\
\hline
\end{tabular}
\caption{Comparison of compute times for a hyperelastic and fractional viscoelastic solutions in $\mathcal{C}$\textbf{heart} with (left) changes in $\Delta_t$ ($ N = 9$) and (right) changes in $ N $ ($ \Delta_t =$ 1e-3).
All times are reported in seconds and based on simulations run using 16 cores on a 3.8 GHz Titan A399 AMD RYZEN Threadripper 32 core system.
}\label{tab:speedelasticviscoelastic}
\end{table}

\end{document}